\documentclass{article}

 \title{Weil-Petersson Completion of Teichm\"{u}ller Spaces and 
Mapping Class Group Actions }
\author{Sumio Yamada\footnote{Supported in part by NSF DMS0096171/0071862}}

\newtheorem{thm}{Theorem}
\newtheorem{lem}{Lemma}
\newtheorem{prop}{Proposition}

\newtheorem{defi}{Definition}
\newtheorem{cor}{Corollary}

\newenvironment{rmk}{{\bf Remark}}{\\}
\newenvironment{pf}{\smallskip {\bf Proof} \smallskip }{\hfill Q.E.D. \smallskip\\}

\newcommand{\teich}{Teichm\"{u}ller}
\newcommand{\weil}{Weil-Petersson}

\newcommand{\T}{\ensuremath{\cal T}}
\newcommand{\Tbar}{\ensuremath{\overline{\cal T}}}
\newcommand{\Tbdry}{\ensuremath{\partial {\cal T}}}
\newcommand{\g}{\ensuremath{\gamma}}
\newcommand{\del}{\partial}
\newcommand{\map}{\ensuremath{{\rm Map}(\Sigma)}}
\newcommand{\e}{\ensuremath{\varepsilon}}

\begin{document}

\maketitle

\begin{abstract}
Given a surface of higher genus, we will look at the 
Weil-Petersson completion of the Teichm\"{u}ller
space of the surface, and will study the isometric action of
the mapping class group on it.   
The main observation is that the geometric
characteristics of the setting bear strong similarities
to the ones in semi-simple Lie group actions on 
noncompact symmetric spaces.
\end{abstract} 

\section{Introduction}


  It is well known~\cite{Wo3} that the Weil-Petersson metric is not
complete on the Teichm\"{u}ller space over a closed  
surface of higher genus.  
When a Weil-Petersson geodesic cannot be further extended, a non-trivial
closed geodesic shrinks in length (with respect to the hyperbolic 
metric) to zero, thus developing a node. 
Take the Weil-Petersson completion $\overline{\cal T}$ of the 
Teichm\"{u}ller space $\cal T$.  It was shown by Masur~\cite{Ma} that
the Weil-Petersson metric extends to $\overline{\cal T}$. 
In this paper, we show that 
the space $(\overline{\cal T}, d)$ is an NPC (or CAT(0)) space in the sense
of Toponogov~\cite {KS1}, even though the distance function $d$  induced by the Weil-Petersson
metric is no longer smooth (with respect to geometric quantities such as
the hyperbolic length of closed geodesics.) 
By construction, the mapping class group (\teich\ modular group)
acts isometrically on the \teich\ space $\T$. One can extend the
isometric action of the mapping class group to the completion $\Tbar$. 
It will be noted that the geometry of $\Tbar$ is closely related to
the isometric actions of various subgroups of the mapping class group. 
Although $\Tbar$ is no longer a manifold, it still has many geometric 
characteristics
shared with the so called Cartan-Hadamard manifolds; complete simply-connected
manifolds with nonpositive sectional curvature.  The similarities
with the action of semi-simple Lie group $G$ on the symmetric space $G/K$
will be noted. To be more specific, the aim of this paper is to rewrite
the paper~\cite{B} of Lipman Bers' where he characterizes, after 
Thurston~\cite{Th}, the elements 
of mapping class group in terms of their translation distances with respect to
the \teich\ metric, only to replace the \teich\ metric by the \weil\ metric.

This paper is motivated to provide a geometric approach to the subject of 
super/strong rigidity where lattices of Lie groups
are represented in the mapping class group of a surface.  As in the 
papers of Corlette~\cite{C}, Gromov-Schoen~\cite{GS}, the rigidity questions can be
transcribed into the study of equivariant harmonic maps into
the NPC space on which the isometry group acts. In this approach,
the negative curvature condition is cruicial to controling analytic
properties of the harmonic maps.  In the
case of strong rigidity, the representation arises as the monodromy
of some fibration where the fiber is the Riemann surfaces of varying
conformal structures.   The monodromy is created by existence of singular
surfaces/fibers, or equivalently vanishing cycles.  
It should be noted that the  
super rigidity of lattices of rank two and higher in mapping
class groups have been studied recently by Farb and Masur~\cite{FM}  via 
a group
theoretic approach. 

Also it should be pointed out that there has been much work done on so-called 
argumented \teich\ space, 
and its mapping class group action on it (see~\cite{Ab} for example). One should note
that the \weil\ completion of a \teich\ space can be identified with the augmented
\teich\ space set-theoretically.    

The author wishes to thank G. Tian for originally suggesting to look 
at the geometry
of moduli space behind the monodromy of Lefschitz fibration, which motivated
this investigation.  He also wishes to thank H. Masur and 
M. Wolf who have offered numerous suggestions and comments
in the course of completing the paper.  G. Daskalopoulos pointed
out a mistake in a previous version of the paper, for which the author is 
grateful.     
And he likes to record 
with gratitude that crucial 
insights in the proof of  \textbf{Theorem 8} were provided by J. Brock 
and also by S. Kerckhoff.        
He thanks the refree for the careful reading of this paper.
Finally the author thanks for R. Schoen's continual support and
encouragement on the project.

\section{Background}


  Let $\Sigma^2$ be a closed (compact and without boundary) surface of genus $g$ with $g > 1$.  Denote the set of all 
smooth Riemannian metrics on $\Sigma$ by ${\cal M}$.  Denote the set of
all hyperbolic metrics on $\Sigma$ by ${\cal M}_{-1}$.  Note that by the
uniformization theorem, ${\cal M}_{-1}$ can be identified with the
set of all conformal structures on $\Sigma^2$.  Let ${\cal D}$ be the group
of smooth orientation-preserving diffeomorphisms of $\Sigma$ , and ${\cal D}_0$ the subgroup of diffeomorphisms
homotopic to the identity map from a fixed Riemann surface $\tilde{\Sigma}$
(this gives markings to all the points in ${\cal M}_{-1}$.)

Define the \teich\  space ${\cal T}_g$ of $\Sigma$ to be 
\[ {\cal T}_g =  {\cal M}_{-1}/{\cal D}_0. \]

Define the moduli space $M_g$ of $\Sigma$ to be 
\[ M_g = {\cal M}_{-1}/{\cal D}. \]

The discrete group ${\cal D}/{\cal D}_0$ is called the mapping class group,
or the \teich\ modular group. which we will denote by $\map$.  

The space $\cal{M}$ of all Riemannian metrics  has a natural $L^2$-metric
defined by \[<h, k>_{L^2(G)} = \int_{N} <h(x), k(x)>_{G(x)} d\mu_G(x)\]
where $h$ and $k$ are symmetric $(0,2)$-tensors, which belong to $T_G \cal{M}$. Knowing that 
${\cal M}_{-1}$ is smoothly imbedded in ${\cal M}$ with the induced $L^2$-metric,
and also that ${\cal M}_{-1} \rightarrow {{\cal M}_{-1}}/{{\cal D}_0}$ is a 
Riemannian submersion (see~\cite{FT}), it makes sense to restrict the $L^2$-metric defined on ${\cal M}$
to ${{\cal M}_{-1}}/{{\cal D}_0}$.  
Thus the Teichm\"{u}ller space has a $L^2$-inner product
structure, and it is called \textit{Weil-Petersson} metric.  It should be noted that
the \weil\ cometric was introduced (Ahlfors~\cite{A})
as an $L^2$ pairing of two cotangent vectors, or equivalently two 
holomorphic quadratic differentials on the surface.  It was then identified
with the $L^2$ metric defined as above by  Fischer and Tromba~\cite{FT}.  
Recall the 
standard geometric fact~\cite{GDH} that any \weil\ geodesic in $\T$ can be lifted
horizontally once the initial point of the lift is specified, and the
lift is then itself a geodesic in ${\cal M}_{-1}$ with respect to
the $L^2$ metric.  In what follows, we will not distinguish a \weil\
geodesic in $\T$ and its horizontal lift in ${\cal M}_{-1}$ unless it is necessary.

With respect to this metric, the Teichm\"{u}ller space
$\mathcal{T}$ has non-positive 
sectional curvature (see  Tromba~\cite{To}, or Wolpert~\cite{Wo2})
and though the metric is  incomplete (Wolpert~\cite{Wo3})  ---not every
 Weil-Petersson geodesic
can be extended indefinitely--- $\mathcal{T}$ is still geodesically convex,
that is, every pair of points can be joined by a unique length minimizing geodsic (Wolpert~\cite{Wo1}.)
It is also known that the space is simply connected, diffeomorphic to the
$6g-6$ dimensional Euclidean ball, where $g$ is the genus of the surface $\Sigma$ (see~\cite{To} for references.)
       
We will first show that the incompleteness is always caused by 
pinching of (at least) one neck of the Riemann surface.  Since the proof 
(as presented in ~\cite{To}) is short and elementary, we will include it here.  

\begin{prop}
Suppose that $\sigma: [0, T) \rightarrow \T$, where $T < + \infty$ 
is a \weil\ geodesic, which cannot be extended beyond $T$. Then for any sequence 
$\{ t_n \}$ with $\lim t_n = T$, the hyperbolic length of the shortest
closed geodesic(s) on $(\Sigma, \sigma(t_n))$ converges to zero. 
\end{prop}

\begin{pf}
Suppose not.  Then there is some lower bound $\e$ for the length of all
closed geodesics in $\Sigma$ on $\sigma ( [0, T) )$.  Then the compactness theorem of Mumford
and Mahler says that there exists a subsequence of $\{ t_n \}$, which we
 donote by $\{ t_n \}$ again, and a sequence of diffeomorphisms $\{ \phi_n \}$
of $\Sigma$ such that $\phi^*_n  \sigma (t_n)$ converges to a hyperbolic
metric $G$.  Note $\phi^*_n \sigma$ is a horizontal lift of a \weil\ goedesic defined on 
$( 0, T]$ for each $n$.  (Here we are using the fact that ${\cal M}_{-1}
\rightarrow \T$ is a Riemannian submersion.)

In the meantime, the existence theorem of ordinary differential equation says 
that given $G$ in the space ${\cal M}_{-1}$ of hyperbolic metrics, there
exist an open neighborhood $U$ of $G$ and  $\delta >0$ such that any geodesic 
with an initial point $G'$ in $U$ is defined on $(-\delta, \delta )$.  

Choose $n$ sufficiently large so that $\phi^*_n \sigma (t_n)$ is in $U$,
and $T-t_n < \delta/2$.  Then the geodesic $\phi^*_n \sigma(t)$
can be extended to the interval $(t_n - \delta, t_n + \delta)$, which
is a contradiction since $T < t_n + \delta$.  
\end{pf}

\begin{defi}
Let $\Tbar$ be the \weil\ completion of the \teich\ space of a Riemann
surface of genus greater than one.  Denote by $\Tbdry$ the 
frontier set $\Tbar \backslash \T$.  
\end{defi}

The preceeding proposition states that every point in $\Tbdry$ represent
a nodal surface, that is, a surface with a node or equivalently a pinched
neck.  H. Masur has shown in~\cite{Ma} that $\Tbdry$ consists of a union
of \teich\ spaces of topologically reduced Riemann surfaces, created by
neck pinching as the conformal structure degenerates toward the frontier
points.  Masur also showed that the \weil\ metric tensor of $\T$ restricted
to the directions tangent to the frontier set $\Tbdry$, spanned by
the holomorpic quadratic differentials with poles of order
one or less over the pinching neck, converges to 
the \weil\ metric tensor of the \teich\ space of the topologically
reduced Riemann surface.  In this sense the \weil\ metric extends to 
$\Tbar$.  The \weil\ metric tensor evaluated in the directions spanned by
holomorphic quadratic differentials with order two poles over the
pinching neck, blows up 
at various rates (also in~\cite{Ma}), which we will carefully analyze in the 
following section.  

Lastly in this section we prove the following theorem, which was
first proved by the author by a different arguement.  It was pointed
out later  by M. Wolf and H.Masur that the statement can be obtained by a direct application of a result (Corollary II 3.11) in the book~\cite{BH}
by Birdson and Heafliger. 

\begin{thm}
The \weil\ completed \teich\ space $\Tbar$ is an NPC space (or equivalently
a CAT(0) space.)
\end{thm}
 
\begin{rmk} 
NPC stands for ``non-positively curved'' as defined in~\cite{KS1}.
It is a length space $(X, d)$, in which any pair of points $p$ and $q$
can be connected by a rectifiable curve whose length realizes
the distance $d(p, q)$, and in which any triangle satisfies 
the length comaprison in the sense of Toponogov with a comparison
triangle in ${\bf R}^2$.
\end{rmk}

\begin{pf}The result (Corrolary 3.11) cited in~\cite{BH} says that the metric completion
of an NPC space is an NPC space.  The \teich\ space equipped with the
\weil\ metric is an NPC space, since it is simply connected, non-positively
curved, geodesically convex, open manifold as described above. Hence 
it follows that its \weil\ metric completion $\Tbar$ is an NPC space.
\end{pf}

\section{Singular Behavior of Weil-Petersson Metric}


We will first consider the case where $P$ in  $\Tbdry$
represents a Riemann surface ${\Sigma}_0$ with one node. It belongs to
a copy of a \teich\ space  ${\T}_{c_1}$ of a  topological surface with a 
node (or equivalently a surface with two punctures.) Suppose  that
this $\Sigma_0$ is obtained by pinching a closed geodesic $c_1$ of 
a non-singular
surface $\Sigma$ (i.e. without nodes) to a point.
Now introduce a complex
coordinate system, as demonstrated in~\cite{Ma}, 
$t = (t_1,..., t_{3g-3})$ where $g$ is the genus of the
non-singular surface $\Sigma$ such that the origin $0$ is $\Sigma_0$,
where $t_2, ... t_{3g-3}$ parametrize the \teich\ space ${\T}_{c_1}$ and 
$t_1$ is induced by the local coordinates near the node $N$ as follows. 

At the node $N$, $\Sigma_0$ has a neighborhood 
isomorphic to $\{ |z| < 1, |w| < 1, zw_1 = 0 \}$ in ${\bf C}^2$.  
Remove two discs $\{ z : 0 < |z| \leq |t_1| \}$ and 
$\{ w : 0 < |w| \leq |t_1| \}$
from $\Sigma_0$, and then identify $z$ with $t_1/w$. We denote by
$\Sigma_t$ the Riemann surface thus obtained.    Given the complex 
structure of $\Sigma_t$, we will assume that $\Sigma_t$ is uniformized, 
that is,
equipped with the hyperbolic metric. As $|t_1| \rightarrow 0$,
the surface $\Sigma$ develops a node $N$.     

Observe that by a pinching a closed geodesic to a point, one can have two
topologically distinct pictures depending on whether $[c_1]$ is
homologically nontrivial or not.  One is when the resulting surface $\Sigma_0$
has one path-connected component, with genus $g-1$ and with two
punctures.  The other is that the surface $\Sigma_0$ consists of two
disconnected surfaces, of genus $g_1$ and $g_2$ with $g_1 + g_2 = g$
and each surface has one puncture. 

In the first case, the frontier component ${\T}_{c_1}$ is the \teich\ space
of surfaces of genus $g-1$ with two punctures.  The complex dimension of 
${\T}_{c_1}$
then is $3[(g-1)-1] +2 = 3g-3 -1$, where the extra two real dimensions
is due to the freedom to choose the positioning of the two punctures.

In the second case, ${\T}_{c_1}$ is a product space of two \teich\ spaces 
${\T}_{c_1}^1$
and ${\T}_{c_1}^2$, where ${\T}_{c_1}^i$ represents the set of Riemann surfaces of genus
$g_i$ with one puncture.  Then the dimension of the product space is
\[
[3(g_1-1)+1] + [3(g_2-1) + 1] = 3(g_1 + g_2 -1)-3+2=3g-3-1.
\]

Hence in either case the dimension of the frontier \teich\ space 
${\T}_1$ is of complex codimension one.  Similarly when $\Sigma_0$ has $n$ nodes, the frontier 
component that parametrizes the nodal surfaces is of complex codimension $n$.

H. Masur~\cite{Ma} showed that the \weil\ metric 
tensor blows up as $|t_1|
\rightarrow 0.$  In particular, he showed that 
\[
0 < \liminf_{t = (t_1, ... t_{3g-3}) \rightarrow 0}  
|t_1|^2 (-\log |t_1|)^3 G_{1 \overline{1}} 
< \limsup_{t = (t_1, ... t_{3g-3}) \rightarrow 0}  
|t_1|^2(-\log |t_1|)^3 G_{1 \overline{1}} < C
\]
where $t = 0 \in {\bf C}^{3g-3}$ represent the surface with the node $P$.
 
We will refine Masur's result and show the following.

\begin{prop} As $|t_1|$ goes down to zero, that is, as a node develops,
one has the following description of the blowing up of the \weil\
metric component.
\[
|G_{1 \overline{1}}(t)| = \frac{1}{\Big( C + O( (- \log |t_1|)^{-2}) \Big)
|t_1|^{2} (- \log |t_1|)^{3}   }.
\]
for some constant $C >0$.   
\end{prop}

\begin{pf}
We will introduce a one-parameter family of hyperbolic surfaces
which models the development of the node as $|t_1|$ goes to zero.

Denote by $A_{|t_1|}$ the annulus $\{ z : |t_1| < |z| < 1 \}$
in ${\bf C}$. One can uniformize the annulus by assigning the following 
conformal factor to the conformal structure of the annulus.
\[
ds^2_{|t_1|} = \Big( \frac{\pi}{\log |t_1|} \csc \frac{\pi \log |z|}{\log |t_1|} \Big| \frac{dz}{z} \Big| \Big)^2.
\]
As the neck pinches ($|t_1| \rightarrow 0$), for each fixed $z$,
the above converges to
the hyperbolic metric on two copies of the punctured disc $\{ 0< |z|  < 1\}$;
\[
ds^2_0 = \Big( \frac{|dz|}{|z| \log |z|} \Big)^2,
\]
which models the standard hyperbolic cusp.

Let $ds^2_{t_1}$ be the hyperbolic metric of $\Sigma_t$ restricted to
the annulus region $A_{|t_1|}$ with respect to the
coordinates defined near the closed geodesic $\{ |z| = \sqrt{|t|} \}$.
Now we quote the following result of Wolf-Wolpert~\cite{WW}.
\[
\| ds^2_{t_1} - ds^2_{|t_1|} \|_{C^2} = O( (- \log |t_1|)^{-2})
\]
on the smaller annulus $A_{|t_1|}^{\delta} =
\{ |t| \leq |z| \leq 1-\delta \}$ for a given
$\delta > 0$.   

Indeed 
 $ds^2_{|t_1|}$ has the following expansion
\[
\begin{array}{lll}
ds^2_{|t_1|} & = & \Big( \frac{|dz|}{|z| \log |z|} \Big)^2 
(\Theta \csc \Theta)^2 \\
             & = & ds^2_0 \Big( 1 + \frac{1}{3} \Theta^2 + \frac{1}{15} \Theta^4+ ... \Big).
\end{array}
\]
where 
\[
\Theta = \pi \log |z| / \log |t_1|,
\]
and for each 
$z$ with $|z| < 1-\delta$, we have
$\Theta(z) \rightarrow 0$ as $|t_1| \rightarrow 0$.            

Now it follows from the Wolf-Wolpert estimate  that
\[
ds^2_{t_1} =  (1 + O( (- \log |t_1|)^{-2} ) ) \frac{\Theta^2}{\sin^2 \Theta } ds^2_0
\]
as $|t_1|$ goes down to zero.  This in turn implies that
$ds^2_t$ converges pointwise to $ds^2_0$ as $|t_1|$ goes down
to zero.  

Now we will follow very closely the method of Masur's to compute the
asymptotics of the \weil\ metric.  Let $\phi_1$ be the cotangent vector
dual to the tangent vector $\del_{t_1}$ against the \weil\ pairing.  Then 
we have the following
expression for $\phi_1$
\[
\phi_1 (z) = - \frac{t_1}{\pi} a(z) \Big( \frac{dz}{z} \Big)^2
\]
near the pinching neck, where $a(z) = 1 + O(|z|)$ for $|t_1| >0$
(see~\cite{Wlf}.)
 
We also need the following estimates~\cite{Ma}.
For any $\delta > 0$, there exists a constant $C > 0$
independent of $t$ such that
\[ 
\frac{1}{C} r^2 (- \log r)^2 \leq  1/\rho^2  (z) \leq 
C r^2 (- \log r)^2
\]
for $\sqrt{|t|} \leq |z|=r \leq 1-\delta$, where 
$\rho^2 (z) = \rho^2 (z, t)$ is the conformal factor as above  
of the complex coordinate $z$ for the 
hyperbolic metric $ds^2_{t}= \rho^2 dz \otimes d \overline{z}$.  

Then take the \weil\ pairing of the cotangent vector $\phi_1$ with itself,
over the modified annulus $A_{|t_1|}^{\delta} =
\{ |t| \leq |z| \leq 1-\delta \}$.
\[
\int_{A_{|t_1|}^{\delta}} \frac{ | \phi_1 |^2 (z) }{ \rho^2 (z)} |z| dr d\theta.
\]

We first note that this quantity is convergent as $| t_1 |$ goes down
to zero.  To see this, first note that as Masur shows, the quantity is 
bounded above and below by $C_1 |t_1|^2 (- \log |t_1|)^3)$
and $C_2 |t_1|^2 (- \log |t_1|)^3)$, with $C_1 > C_2 > 0$.

We claim that the \weil\ pairing above is described by
$ \{ C + O( (- \log |t_1|)^{-2} ) \} |t_1|^2 (-\log |t_1|)^3 $
as $|t_1|$ goes down to zero for some positive constant $C$.

Substituting the expansion of the conformal factor $\rho^2 (z, t_1)$, 
we obtain
\[
\begin{array}{lll}
 & & \int_{A_{|t_1|}^{\delta}} 
\frac{| \phi_1 |^2 (z)}{ \rho^2 (z)} r dr d\theta\\ 
   & = & \{ 1 + O( (- \log |t_1|)^{-2} ) \} |t_1|^2 \int_{A_{|t_1|}^{\delta}} 
   \frac{(1 + O(r))}{r^4} r^2 (-\log r)^2 
   \frac{\sin^2 \Theta^2}{\Theta^2} dx dy
\\
   & = & \{ 1 + O( (- \log |t_1|)^{-2} ) \} |t_1|^2 \int_{A_{|t_1|}^{\delta}} 
   \frac{(1 + O(r))}{r^4} r^2 (-\log r)^2 
   \frac{ \sin^2 \Big( \frac{\pi (- \log r)}{(- \log |t_1|)} \Big)}
   {\Big( \frac{\pi (- \log r)}{(- \log |t_1|)} \Big)^2} dx dy \\
   & = &  \{ 1 + O( (- \log |t_1|)^{-2} ) \} |t_1|^2 \frac{(- \log |t_1|)^{2}}
   {\pi^2}  \int_0^{2 \pi} \int_{|t_1|}^{1 - \delta} \frac{\sin^2 \Big( 
   \frac{\pi (- \log r)}{(- \log |t_1|)} \Big)}{r} d r d \theta_1  
\end{array}
\]
Define $s = \sin \Big( \frac{\pi (- \log r)}{(- \log |t_1|)} \Big)$ and
change the variable.  Then the above is equal to
\[
\{ 1 + O( (- \log |t_1|)^{-2} ) \} |t_1|^2 \frac{(- \log |t_1|)^3}
{\pi^3} \int_0^{2 \pi} \Big[ \int_0^1 \frac{s^2}{\sqrt{1-s^2}} ds 
+ \int_0^{\sin \Big( \frac{\pi (- \log (1- \delta))}{(- \log |t_1|)} \Big)} 
\frac{s^2}{\sqrt{1-s^2}} ds \Big] d \theta
\]
\[
= \{ C + O( (- \log |t_1|)^{-2} ) \} |t_1|^2 (- \log |t_1|)^3
\]
for some positive number $C$. Note here that the number $C$ above 
does not depend on $t_i$'s for it is determined by the value of
the integral $\int_0^1 \frac{s^2}{\sqrt{1-s^2}} ds$. 

As described in Masur's paper, it is known that on any compact set $K$
we have 
\[
\int_K  \frac{| \phi_1 |^2 (z)}{ \rho^2 (z)} dxdy = O(|t_1|^2).
\]
which vanishes faster than the integral of the same integrand over the annulus.
Therefore we have 
\[
\begin{array}{lll}
G^{1 \overline{1}}(t_1) & = & \int_{{\Sigma}_{|t_1|}}
\frac{| \phi_1 |^2 (z)}{ \rho^2 (z)} dxdy \\
 & = & \{ C + O( (- \log |t_1|)^{-2} ) \} |t_1|^2 (- \log |t_1|)^3
\end{array}
\]
The statement of the proposition follows by inverting the matrix 
$G^{i \overline{j}}$ with $j \neq 1$ 
as in the argument given by Masur~\cite{Ma}. Other diagonal terms 
are bounded away from zero.
One point which needs to be addressed in inverting the matrix is 
the fact that all the off-diagonal terms $G^{1 \overline{i}}$ 
vanish at faster rates
than the diagonal terms, and hence the matrix behaves as if it were diagonal.  
\end{pf}

Let $z$ be the local
coordinate near the node employed in the previous argument. Recall
that this particular
coordinate was chosen so that the pinching of the neck is
closely approximated by the hyperbolic cylinder
with the hyperbolic metric $ds_{|t_1|}^2$ defined above.
In particular, it is shown by Wolpert~\cite{Wo4}
that the hyperbolic length $\lambda_1$ of the closed geodesic around the
pinching neck is given by
\[
\lambda_1 (t_1) = \frac{2 \pi^2}{- \log |t_1|} + O \Big( \frac{1}{(- \log |t_1|)^4} 
\Big)
\]
for $t_1$ sufficiently small.
Now let $l_1$ be a new coordinate defined by
\[
l_1 = \frac{2 \pi^2}{- \log |t_1|}
\]

We  would like to know the rate at which $G_{i \overline{j}} (t)$
converges to $G_{i \overline{j}} (0, t_2,...,t_{3g-3})$ 
We will first quote a result of Wolf~\cite{Wlf}, which
says that the hyperbolic metric $g_{t_1} = \rho_{t_1}^2 (z) dz \otimes d 
\overline{z} \ \ (t_1 \neq 0)$
on the  non-degenerate surface $\Sigma_{t_1}$ is real analytic in the
lengths $\vec{\lambda} = (\lambda_1, ..., \lambda_n)$ of the closed geodesics around the 
pinching necks, and as $l_1$ goes to zero, it converges to the cuspidal
hyperbolic metric $g_0 = \rho_0^2 (z) dz \otimes d \overline{z}$ on
$\Sigma_0$. 
Moreover when there is only one
neck pinching Wolf~\cite{Wlf} has shown  that 
\[
|\rho_{t_1}^2 (z) - \rho_0^2 (z)| = O((\lambda_1)^2) = O(\{ \log |t_1| \}^{-2})
\]
over $\Sigma$.  We now claim that the \weil\ metric restricted 
to the directions ${\del}_i$ with $i, j >1$ behaves as follows.   
\[
| G_{ij}(t_1, t_2,..., t_{3g-3}) - G_{ij}(0, t_2,..., t_{3g-3}) |  
=  O ( \{ - \log |t_1| \}^{-2} ).
\]
To see this, take two cotangent vectors $d t_i, d t_j$
with $i, j > 1$, each 
identified with a meromorphic quadratric differentials $\phi_1$
and $\phi_j$ respectively, with at most simple poles at $z = 0$
(see~\cite{Ma}.)
$\phi_i$ and  $\phi_j$ have no other poles away from $z=0$.
 
Then the \weil\ cometric tensor $G^{i \overline{j}}$ is given by
\[
\int_{\Sigma_{t_1}} \frac{\phi_i \overline{\phi_j}}{\rho_{t_1}^2} dx dy. 
\]

As before, we consider the region containing the pinching neck
and the rest separately.  Let $A_{|t_1|}$ be the annulus 
$\{ z : |t_1| < |z| < 1 \}$ in $\Sigma_{t_1}$.  Then recall the Wolf-Wolpert
estimate~\cite{WW} quoted above, which implies with respect to the 
complex coordinate $z$
\[
\| \rho_{t_1}^2 - \rho_0^2  \|_{C^2} = O( \{ - \log |t_1| \}^{-2})
\]
on the smaller annulus $A_{|t_1|}^{\delta} =
\{ |t| \leq |z| \leq 1-\delta \}$ for a given
$\delta > 0$.   Hence we have
\[
\begin{array}{lll}
\int_{A_{|t_1|}^{\delta}} \frac{\phi_i \overline{\phi_j}}{\rho_{t_1}^2} dx dy
- \int_{A_{|t_1|}^{\delta}} \frac{{\phi}_i \overline{\phi_j}}{\rho_0^2} dx dy & = & 
\int_{A_{|t_1|}^{\delta}} \phi_i \overline{\phi_j} \Big[ \frac{1}
{[1 + O( \{ - \log |t_1| \}^{-2})] \rho_0^2 (z)} - 
\frac{1}{\rho_0^2 (z)} \Big] dx dy \\
 & = & \int_{A_{|t_1|}^{\delta}} O( \{ - \log |t_1| \}^{-2} )
 \frac{\phi_i \overline{\phi_j}}{\rho_0^2 (z)} dx dy\\
 & = & O( \{ - \log |t_1| \}^{-2} )
 \end{array}
\]
where the term $O( \{ - \log |t_1| \}^{-2} )$ is bounded 
in terms of $\{ - \log |t_1| \}^{-2}$, uniformly in 
$z$ in $A_{|t_1|}^{\delta}$ due to the fact that the Wolf-Wolpert 
estimate is a $C^2$ (in particular $C^0$) 
estimate on $A_{|t_1|}^{\delta}$. The last equality follows
from the fact that the part of the integrand $\phi_i \overline{\phi_j}$
is a term which as $z \rightarrow 0$ can blow up no faster than
the rate of $1/|z|^2$, which in turn implies that the integral
$\int_{A_{|t_1|}^{\delta}} 
\frac{\phi_i \overline{\phi_j}}{\rho_0^2 (z)} dx dy$ is a
term $O(1)$ as $t_1$ goes to zero.

On a compact set $K$ away from the pinching neck, we have
\[
\begin{array}{lll}
\int_K \frac{\phi_i \overline{\phi_j}}{\rho_{t_1}^2} dx dy
- \int_K \frac{\phi_i \overline{\phi_j}}{\rho_0^2} dx dy
& = & \int_K \phi_i \overline{\phi_j} \frac{\rho_0^2 (z) 
- \rho_{t_1}^2 (z)}{\rho_0^2 (z) \rho_{t_1}^2 (z)}
dx dy \\
& = & O( \{ - \log |t_1| \}^{-2} ) \int_K  
\frac{\phi_i \overline{\phi_j}}{\rho_{t_1}^2 (z)}
dx dy \\
 & = & O( \{ - \log |t_1| \}^{-2} )
\end{array}
\]
where the second equality follows from the fact 
shown by Wolf~\cite{Wlf} as already described above, 
that 
$\rho_{t_1}^2$ is real analytic in $\lambda_1 = 2 \pi^2 
\{ - \log |t_1| \}^{-1}
+ O(\{ - \log |t_1| \}^{-4})$ and $\rho_{t_1} (z) 
- \rho_0 (z) = O(\lambda_1^2)$ pointwise on $\Sigma$. The last equality 
follows from the fact that 
the integrand of the previous line is continuous in $z$ over $K$.

Combining those estimates, we see that the difference between 
$G^{i \overline{j}} (t_1, t_2,..., 
t_{3g-3})$ and $G^{i \overline{j}} (0, t_2,..., 
t_{3g-3})$ is a term of $ O(\{ - \log |t_1| \}^{-2})$.

Introduce a new variable $u_1 = \sqrt{l_1}$ 
here.  Then the description
of the \weil\ metric near the frontier point is written down as
\[
\begin{array}{lll}
ds^2 & = & \{ C + O((u_1)^4 \}
d u_1^2 +   \frac{1}{4} \{ C + O((u_1)^4) \} (u_1)^6 d \theta_1^2\\
  & & \ \  \ + \{ \tilde{C} +O((u_1)^4) \} (u_1)^3 
  \times \Big[ \mbox{cross terms of $d u_1$ and $dt_i$'s 
(or $d \overline{t_i}$) } \Big] \\
  & &  \ \ \ \ +  \{ \hat{C} + O((u_1)^4) \} (u_1)^6 \times \Big[ \mbox{cross terms of $d \theta$ and $dt_i$ 
(or $d \overline{t_i}$) } \Big] \\
 & & \ \ \ \ \ \ + \sum_{1 < j \leq 3g-3} \Big( 1 + O( (u_1)^{4} ) \Big) 
     |d t_j|^2 
\end{array}
\]
with $i > 1$. Here we have used the following relations due
to the change of variables.
\[ t_1 = |t_1| e^{i \theta_1}, \ \ \ l_1 = \frac{2 \pi^2}{-\log |t_1|}
\mbox{ and } \ u_1 = \sqrt{l_1} \]
\[ dt_1 = e^{i \theta_1} d |t_1| + i t_1 d \theta_1, \ \ \ 
d \overline{t_1} = e^{-i \theta_1} d |t_1| - i \overline{t_1} d \theta_1 \]
\[ 
\begin{array}{lll}
\Re \Big[ \frac{1}{|t_1|^2 (-\log |t_1|)^3} dt_1 
\otimes d \overline{t_1} \Big] & = & 
\frac{1}{|t_1|^2 (-\log |t_1|)^3} \Big[ (d |t_1|)^2 + |t_1|^2 (d \theta_1)^2
\Big] \\
 & = & \frac{1}{(-\log |t_1|)^3} \Big( \frac{d |t_1|}{|t_1|} \Big)^2 +
 \frac{1}{(-\log |t_1|)^3} (d \theta_1)^2 \\
 & = & \frac{(d l_1)^2}{l_1} + (l_1)^3 (d \theta_1)^2 \\
 & = & 4 (du_1)^2 + (u_1)^6 (d \theta_1)^2
 \end{array}
\]
where $\Re$ denotes the real part of the complex-valued tensor.
\begin{prop}
The \weil\ metric tensor near the frontier ${\T}_{c_1}$ is continuously differentiable 
in $u_1, \theta_1$ and $t_i$'s.  
\end{prop}

\begin{pf}
We will show that the 
\weil\ cometric tensor $G^{ij}$ is continuously differentiable 
by showing that the differentiations with respect to the 
parameters $(- \log |t_1|)^{-1}$, $\theta_1$, and $t_i$'s
commute with the integration over the hyperbolic surface.
Since $u_1$ is defined to be $\sqrt{ 2 {\pi}^2 (- \log |t_1|)^{-1}}$, 
the statement of the proposition then follows.

In \cite{WW} Wolf and Wolpert showed that the hyperbolic metric
is sector-real-analytic in $(- \log |t_1|)^{-1}$, $\theta_1$, and $t_i$'s,
that is, for any ray from the origin $t=0$ there is a sector of 
the neighborhood at the origin containing that ray in which the
tensor has a convergent expansion in the variables.  

As for $G^{1 \overline{1}}$ the expression for \weil\ pairing 
of a deformation tensor $\phi$ with
itself over the
annulus region $A_{|t_1|}^{\delta}$;
\[
\begin{array}{lll}
 & & \int_{A_{|t_1|}^{\delta}} 
\frac{| \phi_1 |^2 (z)}{ \rho^2 (z)} |z| d|z| d\theta\\ 
   & = & |t_1|^2 \int_{A_{|t_1|}^{\delta}} \{ 1 + O( (- \log |t_1|)^{-2} ) \} 
   \frac{(1 + O(|z|))}{|z|^4} |z|^2 (-\log |z|)^2 
   \frac{ \sin^2 \Big( \frac{\pi (- \log |z|)}{(- \log |t_1|)} \Big)}
   {\Big( \frac{\pi (- \log |z|)}{(- \log |t_1|)} \Big)^2} dx dy \\
   & = &  |t_1|^2 \frac{(- \log |t_1|)^{2}}
   {\pi^2}  \int_0^{2 \pi} \int_{|t_i|}^{1 - \delta} 
   \{ 1 + O( (- \log |t_1|)^{-2} ) \} \frac{\sin^2 \Big( 
   \frac{\pi (- \log |z|)}{(- \log |t_1|)} \Big)}{|z|} d |z| d \theta 
\end{array}
\]
where the term $\{ 1 + O( (- \log |t_1|)^{-2} ) \}$ is sector-real-analytic
in $(- \log |t_1|)^{-1}$, $\theta_1$ and $t_i$'s.  Denote the integrand 
of the last line above by $f(z, (- \log |t_1|)^{-1}, \theta_1, t_i)$.
Formally we can differentiate the integral
\[
\int_0^{2 \pi} \int_{|t_i|}^{1 - \delta} 
   f(z, (- \log |t_1|)^{-1}, \theta_1, t_i) d|z| d \theta_1 
\]
with respect to $(- \log |t_1|)^{-1}$ and obtain
\[
\frac{\del }{\del (- \log |t_1|)^{-1}} \int f d|z| d \theta_1 = \int_0^{2 \pi} 
\int_{|t_i|}^{1 - \delta} 
\frac{\del f}{\del (- \log |t_1|)^{-1}} d|z| d \theta_1 - 
\int_0^{2 \pi} f \Big|_{|z|=|t_1|} d \theta_1 
\]
This expression is justified as follows.
The first term exists since the difference quotient
\[
\frac{f(z, (- \log |t_1|)^{-1}+\e) - f(z, (- \log |t_1|)^{-1})}{\e}
\]
is uniformly bounded as $\e$ goes down to zero, so that one
can apply the Lebesgue dominance convergence theorem (note here 
we need to extend the 
domain of $f$ in $z$ suitably to take the difference.)  The second term 
is zero since when $ |z| = |t_1|$, we have $f(z) = 0$. 

Similarly we can differentiate the same term with respect to $\theta_1$
and $t_i$'s.  Note that the only dependence of $f$ on $\theta_1$ and
$t_i$'s comes from the term $\{ 1 + O( (- \log |t_1|)^{-2} ) \}$ 
and from the result of~\cite{WW} we know that this term is differentiable
in those variables.  

On the complement of the pinching neck, the sector-real-analytic dependence
of the hyperbolic metric on all the variables once again 
induces the differentiability of the \weil\ pairing.

The differentiability of $G^{ij}$ can be checked analogously.

\end{pf}

Now we turn our attention to the case where the frontier point $P$ in $\Tbar$
represents a Riemann surface with more than one node.   Let $p > 1$ be 
the number of nodes.  Recall that $p$ is bounded by $3g-3$; the maximal
number of mutually disjoint closed geodesics on $\Sigma$.

Having the convergence of the hyperbolic metrics as a node develops
as studied in the proof of the previous proposition, we can now improve the 
estimates of Masur's and get the following blow-up rates of the 
\weil\ metric tensor.

\begin{prop}
In the neighborhood of $t = (0, 0)$ in ${\mathbf{C}}^{3g-3} = {\bf{C}}^p \times
{\bf{C}}^{3g-3-p}$, where $t_i \ \ (1 \leq i \leq p)$ parametrizes the 
sizes of the pinching necks,  the \weil\ metric is parametrized as follows; as
$t = (t_1, t-2,..., t_{3g-3}) \rightarrow 0$, 
\begin{tabbing}
1) $|G_{i \overline{i }}(t)| = \Big(C + O((- \log |t_1|)^{-2}) \Big) \Big[ 
|t_1|^{2} (- \log |t_1|)^{3} \Big]^{-1}$ for $1 \leq i \leq p$ \\
2) $|G_{i \overline{j}}(t)| = \Big(C  + O((- \log |t_i|)^{-2}) \Big) 
\Big[ |t_i||t_j|(- \log |t_i|)^3 
(- \log |t_j|)^3 \Big]^{-1} $ for $1 \leq i, j, \leq p$ and  $i \neq j$. \\
3) $ \Big| G_{i \overline{j}}(t_1,...,  t_{3g-3}) - G_{i \overline{j}}
(0,..., 0, t_{p+1}, ..., t_{3g-3}) \Big| = O \Big( \sum_{k=1}^p
(- \log |t_k|)^{-2} \Big)$ for $i, j > p$ \\
4) $|G_{i \overline{j}} (t)| = \Big( C + O ( \{ - \log |t_i| \}^{-2}) \Big)  
\Big[ |t_i| (- \log |t_i|)^3 \Big]^{-1} $ for $i \leq p$ and $j > p$.
\end{tabbing}

\end{prop}

In proving the proposition, Masur's proof is modified  
at two technical points: the first is  whenever
there is an integration over a pinching neck the finer convergence of
the hyperbolic metric parametrized by the $t_i$'s is used, and the second  
technical improvement is to use the convergence of the hyperbolic metric
away from the nodes using the estimates by Wolf as quoted above.   

As before we perform the change of variables, this time set
$u_i = \sqrt{l_i} \ \ (1 \leq  i \leq p)$.  Then the \weil\ metric near
the origin $(0,..., 0)$ in ${\mathbf{C}}^{3g-3}$ 
representing the nodal surface $P$ with p nodes has the following 
expression;
\[
\begin{array}{lll}
ds^2 & = & \sum_{i = 1}^{p} \Big(C_i + O( (u_i )^{4} ) \Big)
     d u_i^2  \\
     &   & \ \ \  + \sum_{p< j \leq 3g-3} \Big( 1 + O( (u_j)^{4} \Big) 
     |d t_j|^2  \\
     &  & \ \ \ \ \  + \sum_{1\leq i, \ j \leq p} \Big( \tilde{C_{ij}}+ 
     O( ( u_i )^{4}) + O( ( u_j )^{4}) \Big) (u_i)^3  (u_j)^3
      \times \Big[ \mbox{ cross terms of $du_i$ 
     and $du_j$ } \Big]  \\
     &  & \ \ \ \ \ \ \  + \sum_{k \leq p, \ l>p} 
     \Big( \hat{C_{kl}} +O((u_k)^4) \Big) (u_k)^3 \times 
     \Big[ \mbox{cross terms of $du_k$ and $dt_l$ (or  $d \overline{t_l}$) } 
     \Big]   \\
     &  & \ \ \ \ \ \ \ \ + \sum_{1 \leq i \leq p, \ p< j \leq 3g-3} 
      \Big( \overline{C_i} + O( (u_i )^{4} )  \Big)  (u_i)^6
     \times \Big[ \mbox{ cross terms of $d \theta_i$ and $d t_j$ (or $d 
     \overline{t_j}$)} \Big]  \\
     & & \ \ \ \ \ \ \ \ \ \ \  + \sum_{i = 1}^p \frac{1}{4}\Big( C_i
      + O((u_i)^4) \Big) (u_i)^6 d \theta_i^2 + 
     \sum_{1 \leq i,j \leq p} O((u_i)^6 (u_j)^6) d \theta_i 
     \otimes d \theta_j. 
\end{array}
\]

\section{Geometry of the Frontier Set $\Tbdry$}


We start this section with a theorem which describes how each boundary component is
embedded in $\Tbdry$.

\begin{thm}
Each component of the boundary \teich\ spaces is totally geodesic; that is, given any pair of 
points $p$ and $q$ in a \teich\ space ${\T}_C$ representing a collection of nodal surfaces
${\Sigma}_C$
obtained by pinching a collection $C$ of mutually disjoint simple closed geodesics $c_i$ of
the angular surface $\Sigma$, a length minimizing geodesic connecting $p$ and $q$
are totally contained in ${\T}_C$ and it is unique.  
\end{thm} 

\begin{pf}
Suppose $C = \cup_{i=1}^{|C|} c_i$. 
Let $l_{c_i}(x)$ be the hyperbolic length of the simple closed geodesic
$c_i$ with respect to the hyperbolic metric $x$ on $\Sigma$.  The domain of
the functional $l_{c_i}$ can be continuously extended to $\cup_{A \subset C}
 {\T}_A$ from $\T$ by defining $l_{c_i}|_{{\T}_A} \equiv 0$ if $c_i \in A$.

Define a new functional $L_C: \cup_{A \subset C}
 {\T}_A \rightarrow {\bf R}$ 
by
\[
L_C (x) = \sum_{i=1}^{|C|} l_{c_i} (x).
\]
Note that $L_C|_{{\T}_C} \equiv 0$, hence that $L_C (p) = L_C (q) = 0$.

We now construct a length minimizing geodesic connecting $p$ and $q$.  
Let $\{ p_i \}$ and $\{ q_i \}$ be  Cauchy sequences in $\T$ converging 
to $p$ and $q$ respectively.  Let $\sigma_i (t)$ be the unique
length minimizing \weil\ geodesic connecting $p_i = \sigma_i (0)$ and 
$q_i = \sigma_i (1)$. Note that $\sigma_i$ lies entirely in $\T$ due to 
the geodesic convexity of $\T$~\cite{Wo5}. Then by the strictly
negative sectional curvature of the \weil\ metric on $\T$, we know that
\[
d(\sigma_i (t), \sigma_j (t)) \leq \max \Big( d(p_i, p_j), d(q_i, q_j)  
\Big).
\]
The right hand side of the inequality converges to zero, and hence
it follows that $\sigma_i (t)$ converges to a point in $\Tbar$, which we 
call $\sigma (t)$. 

Consider the composite function $f(t) = L_C (\sigma(t))$.   
We now use the fact that the length functional $l_{c_i}$ is
convex with respect to the \weil\ metric on $\T$ and ${\T}_{c_i}$ (a result of 
S. Wolpert~\cite{Wo5}, see~\cite{Y1} for generalizations.)
Then it follows that $L_C$ is convex on $\T$ and ${\T}_C$, and  hence that
$f(t)$ is convex in $t$.  Suppose that $M = \max f(t) > 0$. Then 
it follows that $f(t) \equiv M >0$ which contradicts
with $f(0) = f(1) = 0$. 

We have so far shown that $f(t) \equiv 0$,
which then implies that $\sigma$ lies in $\overline{{\T}_C}$
To see $\sigma$ lies in ${\T}_C$, note that there exists 
a \weil\ length realizing geodesic $\sigma'$ connecting
$p$ and $q$ lying entirely in the \teich\ space ${\T}_C$  due to the \weil\ geodesic
convexity of the ${\T}_C$.  Since given two points in an NPC/CAT(0) space $X$  a length-realizing geodesic  is unique (here we take $X$ to be 
${\T}_C$ ), we know that $\sigma'$ is {\it the} geodesic $\sigma$ connecting $p$ and $q$.

\end{pf}

We are in a position to present a series of results which illustrate the 
geometry of the frontier set $\Tbdry$.  

\begin{prop}
Let $q$ be a nodal surface, and $p$ a non-degenerate surface.  Then the
\weil\ open geodesic segment connecting $p$ and $q$
lies entirely in the interior \teich\
space $\T$.
\end{prop}

\begin{pf}
We suppose the contrary; that we have a length minimizing geodesic $\sigma$
connecting $p$ and $q$ with a part of $\sigma$ lying in the frontier, and
then moves into the interior \teich\ space to reach $p$.

We will consider the case that $q$ represent a degenerate surface 
with one node with a simple closed geodesic $c$ pinched first.
Choose the coordinates so that the origin corresponds to the point 
$\sigma (0) = r$
at which $\sigma$ leaves the frontier.  Then the geodesic segment
$\overline{pq}$ lies in the frontier and the rest is in $\T$ by the
geodesic convexity of the spaces $\T$ and ${\T}_c$.
Let the $3g-4$ complex dimensional linear space $(t_2, ..., t_{3g-3})$ be 
the geodesic normal coordinate centered at the origin
(hence the geodesic segment $\overline{qr}$ is a ray starting from the
origin.) Let the vertical axis be parametrized by $u = C(- \log |t_1|)^{-1/2}$.
where $C$ is chosen so that the \weil\ metric tensor at the origin restricted
to the plane spanned by ${\del}_x$ and ${\del}_u$. 
is approximated by a model metric $du_1^2 + 1/4 (u_1)^6 d \theta_1^2 + 
\sum_{i = 2}^{3g-3} dt_i \otimes d \overline{t_i}$ 
defined on the upper-half space $\{ u_1 \geq 0 \}$.

Now consider the geodesic connecting $p$ and $q$ as an one-dimensional
harmonic map $v: [-1, 1] \rightarrow \Tbar$ with the Dirichlet boundary
condition $v(-1) = p$ and $v(1) = q$. 
The next step is the following lemma.

\begin{lem} 
The pull-back $v^* u_1 = u_1 (v)$ of the coordinate function 
$u_1 = C/\sqrt{- \log |t_1|}$
by the harmonic map $v$ satisfies the differential inequality 
\[
\frac{d^2}{dt^2} (v^* u_1)(t) \leq C_1 (v^* u_1)(t)
\]
on $(-1, 1)$ distributionally for sufficiently small values of 
$v^* u_1$ for some
constant $C_1$.
\end{lem}

Once we have this estimate, we can proceed to prove the
theorem as follows.  The next lemma is the one dimensional version of 
Harnack inequality for $W^{1, 2}$ functions (functions whose derivatives 
are in $L^2$.)   

\begin{lem}Suppose that non-negative $W^{1,2}$ function $u$, defined
on $B_{2R_0} := \{ -2 R_0 < x < 2 R_0 \}$, satisfies the  
differential inequality $u'' < C_1 u$ weakly. Then $u$ satisfies    
\[
\sup_{B_{R_0}} u < C_2 \inf_{B_{R_0}} u  
\]
for some constant $C_2>0$ independent of $u$ and $R_0$.
\end{lem}

Note that  $\inf_{B_{R_0}(0)} (v^* u_1)=0$ since $v(0) = r$
denote a nodal surface, and the inequality then
implies $v^* u_1 \equiv 0$ on $B_{R_0}(0)$ which is a contradiction
to the initial supposition that $(v^* u_1)(t) >0$ for $t>0$.

\begin{pf}[of the Harnack-type inequality]
Modify $u$ so that
\[
w_{\varepsilon}(x) = u(x)  + \varepsilon.
\]
We have $\varepsilon > 0$ fixed, and later let it go to zero.
In the meantime, we suppress the dependence on $\varepsilon$;
$w = w_{\varepsilon}$. Note that $w$ satisfies the same 
differential inequality $w'' \leq C_1 w$ as $f$.
Choose $ \zeta $ to be a compactly supported cut-off function 
on $B_{2 R_0}$ such that 
$\zeta \equiv 1$ on $B_{R_0}$, and $|\zeta'(x)| < 2/ R_0$ on $B_{2 R_0} \backslash B_{R_0}$.
\[
\begin{array}{lll}
\log \frac{\sup_{B_{R_0}} w}{\inf_{B_{R_0}} w} & = & \log ( \sup_{B_{R_0}} w ) - \log ( \inf_{B_{R_0}} w)\\ 
  & = &  \sup_{B_{R_0}}  \log w  - \inf_{B_{R_0}} \log w\\
  & \leq & \int_{B_{R_0}} | (\log w)'(x) | dx \\
  & \leq & \Big( \int_{B_{R_0}} | (\log w)'(x) |^2 dx \Big)^{1/2} \sqrt{2 R_0} \\
  & \leq & \Big( \int_{B_{2 R_0}} \zeta^2 (x) | (\log w)'(x) |^2 dx \Big)^{1/2}
 \sqrt{2 R_0} 
\end{array}
\]
We now claim that 
\[
\int_{B_{2 R_0}} \zeta^2 (x) | (\log w)'(x) |^2 dx < C_3 / R_0
\]
for any $0 < R_0 < M$ and $C_3 = C_3 (M)$.  Once we have this estimate, we 
combine the two inequalities above
we get for $0 < R_0 < M$,
\[
 \log \frac{\sup_{B_{R_0}} w}{\inf_{B_{R_0}} w} < C_2.
\]
Finally we let $\varepsilon$ goes to zero and the statement of the lemma follows. 

To see the claim hold true, note that 
\[
\begin{array} {lll}
\int_{B_{2 R_0}} C_1 w (\zeta^2 w^{-1}) dx  & \geq  & \int_{B_{2 R_0}} w'' (\zeta^2 w^{-1}) dx \\
 & = & - \int_{B_{2 R_0}} [ -w^{-2} (w')^2 \zeta^2 + 2w' w^{-1} \zeta {\zeta}' ] dx \\
  & = & \int_{B_{2 R_0}} \Big[ (\log w)' \Big]^2 \zeta^2 - 2 (\log w)' \zeta \zeta' dx
\end{array}
\]
The second integral should be regarded as $w''$ being integrated against a test function
$(\zeta^2 w^{-1})$, and the first inequality is due to the hypothesis $u'' \leq C_1 u$. By
reorganizing the terms, we get 
\[
\begin{array}{lll}
\int_{B_{2 R_0}} \Big[ (\log w)' \Big]^2 \zeta^2 dx & \leq & \int_{B_{2 R_0}} 2 (\log w)' \zeta \zeta' +  C_1 \zeta^2  dx \\
  & \leq & \int_{B_{2 R_0}} \frac{1}{2} \Big[ (\log w)' \Big]^2 \zeta^2  + 2 (\zeta')^2 +  C_1 \zeta^2  dx \\
\end{array}
\]
The second inequality is obtained by the inequality $2ab \leq a^2 + b^2$.

Hence it follows that
\[
\int_{B_{2 R_0}} \Big[ (\log w)' \Big]^2 \zeta^2 dx \leq 2 \int_{B_{2 R_0}} 2 (\zeta')^2 +  C_1 \zeta^2  dx < \frac{C_3}{R_0}
\]
since $|\zeta'| < 2/R_0$ and $0 \leq \zeta \leq 1$ on $B_{2R_0}$.
\end{pf}

\begin{pf}[of {\bf Lemma 1}]
The harmonic map equation (or the geodesic equation)
for $v: [-1, 1] \rightarrow \Tbar$
for the first coordinate function $v^1$ is
\[
(v^1)'' + \Gamma^1_{\alpha \beta}(v(t)) (v^{\alpha})' (v^{\beta})'  = 0
\]
where $\Gamma^1_{\alpha \beta}$ is the Christoffel symbol for the
\weil\ metric.  Note that $v^1 (t)$ is nothing but the pulled-back
function $(v^* u_1)(t)$.  (Also $v^2 (t) = (v^* \theta_1)(t)$.) 
Hence to show the inequality, it is
equivalent to showing that the terms
\[
\sum_{\alpha, \beta} \Gamma^1_{\alpha \beta}(v(t)) (v^{\alpha})' (v^{\beta})' 
< C (v^* u_1)(t)
\] 
for some $C > 0$. In fact, we will show 
\[
 \sum_{\alpha, \beta} \Gamma^1_{\alpha \beta}(v(t)) (v^{\alpha})' (v^{\beta})' 
= O (v^* u_1 (t))
\]
which allow us to choose $C > 0$ for sufficiently small $v^* u_1$.

Recall from the previous section that with respect to the coordinate 
system $t = (u_1, \theta_1, t_2,...,t_{3g-3})$
near $t=0$, the \weil\ metric tensor has the following form;
\[
\begin{array}{lll}
G_{11}(t) & = & 1 + O( (u_1)^4 ) \\
G_{12}(t) & = & 0 \\
G_{1j}(t) & = & O( (u_1)^3 ) \ \ (j > 2) \\
G_{22}(t) & = & \frac{1 + O( (u_1)^4 )}{4} (u_1)^6
\\
G_{2j}(t) & = & O((u_1)^6) \ \ (j > 2)\\
G_{ij}(t) & = & (1 + O((u_1)^4)) G_{ij}(0, 0, t_2,...,t_{3g-3}) \ \ 
\mbox{ for $i, j > 2$}
\end{array}
\]
as well as 
\[
G^{11} = O(1), G^{12} = O(1) \mbox{ and } G^{1j} = O((u_1)^3) \ \ 
\mbox{for $j >2$}.
\]

The Christoffel symbols are obtained from the metric tensor by the
following formula (see~\cite{GDH} for example)
\[
\Gamma^1_{ij} = \frac{1}{2} \sum_{l} G^{1l} (G_{il, j} + G_{lj, i} - G_{ij,l}).
\]

Recall here that in {\bf Proposition 3} we showed that the metric tensor is continuously
differentiable with respect to the variables $u_1, \theta_1, t_i's$ and $\overline{t_i}'s$.  
Hence in differentiating the $O((u_1)^k)$ terms appearing the description of $G_{ij}'s$, 
the first derivatives behave as follows; 
\[
\frac{\del O((u_1)^k)}{\del u_1} = O((u_1)^{k-1})
\]
\[
\frac{\del O((u_1)^k)}{\del x} = O((u_1)^k).
\]
where $x$ is any one of the variables other than $u_1$.  This holds because 
$O((u_1)^k) = f(x) (u_1)^k + o((u_1)^k)$ where $f(x)$ is a continuously 
differentiable function {\it independent} of $u_1$ (though it depends on 
other variables $x$.)

Consequently we have
\[
\Gamma^1_{11} = O((u_1)^3), \ \ \Gamma^1_{12} = O((u_1)^4), \ \ \Gamma^1_{22}
= O((u_1)^5)
\]
\[ 
\Gamma^1_{1j} = O((u_1)^2), \ \ \Gamma^1_{2j} = O((u_1)^3) \ \
\mbox{ and } \ \ \Gamma^1_{ij} = O((u_1)^3) \mbox{ for $i, j > 2$}.
\]

Note that the energy density of the harmonic map $v$ is uniformly 
bounded, say by $M^2>0$  (in fact constant, since the geodesic is parametrized
by the arc-length) we know that 
\[ 
|(v^{\alpha})'| < M < \infty
\]
for $\alpha \neq 2$ and for $v^{2} = \theta_1 (v)$ we have 
\[
(u_1)^6 (\theta')^2 < M < \infty.
\]
It is easy to see that $|(\theta_1)'| = O(1)$ as $u_1$ goes to zero,
for if $|(\theta_1)'| = O((u_1)^{-\kappa})$ for some $\kappa >0$,
a comparison map $\tilde{v}$ where ${\tilde{v}}^2 = (\theta_1)(\tilde{v})
= \mbox{constant}$ and otherwise $\tilde{v}$ is defined identical to 
$v$ has less energy, which contradicts $v$ being energy minimizing.

Now with the estimates above note that the term out of the geodesic equation
\[
\sum_{\alpha, \beta} \Gamma^1_{\alpha \beta}(v(t)) (v^{\alpha})' (v^{\beta})' 
\]
is of the order $O((u_1)^2)$.  Hence for $u_1 < \e$ for a sufficiently small
$\e > 0$, we have the differential inequality
\[
\frac{d^2}{dt^2} (v^* u_1)(t) \leq C (v^* u_1)(t)
\] 
for $(v^* u_1)(t) > 0$ where $v$ is a smooth map.  

We will show that the inequality is valid over the
extended region $u_1 \geq 0$ distributionally, that is for any non-negative
smooth compactly supported test function $\phi$ on $[-1, 1]$
\[
\int_{-1}^1 (v^* u_1) \phi'' dt \leq \int_{-1}^1 C (v^* u_1) \phi dt.
\]
Recall that $(v^* u_1)(t) = 0$ for $-1 \leq t \leq 0$ and that
$(v^* u_1)(t) > 0$ for $0 < t \leq 1$. Hence
\[
\begin{array}{lll}
\int_{-1}^1 (v^* u_1)(t) \phi'' dt & = & \int_{-1}^0 (v^* u_1)(t) \phi'' dt
+ \int_0^1  (v^* u_1)(t) \phi'' dt \\
 & = & -\int_0^1  (v^* u_1)'(t) \phi' dt + (v^* u_1) \phi' \Big|_0^1 \\
 & = & \int_0^1 (v^* u_1)'' \phi dt + (v^* u_1)' \phi \Big|_0^1 \\
 & \leq & \int_0^1 C (v^* u_1) \phi dt \\
 & = & \int_{-1}^1 C (v^* u_1) \phi dt. 
\end{array}
\]
where we have used the facts that $(v^* u_1)(0)=0, \phi'(1) = \phi (1) = 0$,
$(v^* u_1)' (0) = 0$ and that $v$ is Lipschitz continuous.    
The last equality $(v^* u_1)' (0) = 0$
is due to the following argument taken out of~\cite{GS}. 
By taking a sequence of scalings of both the domain metric
and the target distance function at $0$ and $v(0)$ respectively, 
one obtains the homogeneous map $v_* (s)$
from ${\bf R}$ to the tangent cone of $\Tbar$ at $v(0)$, which is by itself
harmonic since dilations preserves the harmonicity of the map.  
The tangent cone is isomorphic to ${\bf R}^{6g-8} \times {\bf R}^+$.
Now since $t \in {\bf R}$, $v_*(t) \in {\bf R}^{6g-8} \times \{ 0 \}$,
it follows that  $(v^* u_1)' (0) = 0$.  This is because if  
$(v^* u_1)' (0) > 0$, then after taking the limit of the scaling, 
the resulting map $v_*$ would not be harmonic, for the image of ${\bf R}$ 
will have a corner at the origin of the tangent cone, (namely a jump
discontinuity in its first derivative) a contradiction to the
fact that on ${\bf R}$ linear functions are the only harmonic functions.

The fact that $v$ is Lipschitz continuous was also shown in~\cite{GS}
by showing that the Bochner's formula holds distributionally, which
then implies the DiGiorgi-Nash-Moser type estimate applied to the energy
density of the geodesic map $v$, showing that the modulus of 
continuity of $v$ is bounded.

This proves the lemma.
\end{pf}

To complete the proof of the proposition, we need to consider the case when 
the point $q$ represents a nodal surface $\Sigma_C$ where $C = \cup_i^N c_i$
is a collection of more than one simple closed mutually disjoint
geodesics $c_i$'s. Recall that $p$ is a point in the interior \teich\ space
$\T$.  Suppose now that the geodesic $\sigma$ starting at $q = \sigma (-1)$ 
travels within ${\T}_C$ till leaves it at $r = \sigma (0)$ in ${T}_C$.
As before, denote by $v$ the map $v : [-1, 1] \rightarrow \Tbar$
with the Dirichlet condition $\sigma(-1) = q$ and $\sigma(1) = p$.
Now consider the pull-back function $v^* (\sum_{i=1}^N u_i)
= \sum_{i=1}^N u_i (v(t))$.  Since each function $v^* u_i$
satisfies the weak differential inequality
\[
\frac{d}{dt^2} (v^* u_i)(t) \leq C_i (v^* u_i)(t)
\]
over any open interval $(a, b) \subset [-1, 1]$ for some $C_i >0$. 
By taking the sum over $i$ of the inequalities,
\[
\frac{d^2}{dt^2} \Big[v^* (\sum_{i=1}^N u_i) \Big](t) \leq C 
\Big[ v^* (\sum_{i=1}^N u_i) \Big](t)
\]
for $C = \max C_i$.

Therefore as before we have the following Harnack type inequality
\[
\sup_{(-1/2, 1/2)} \Big[v^* (\sum_{i=1}^N u_i) \Big] \leq C
\inf_{(-1/2, 1/2)} \Big[v^* (\sum_{i=1}^N u_i) \Big]
\]

Given that $\sigma(0)$ represent a surface with multiple nodes, or
equivalently  $[v^* (\sum_{i=1}^N u_i)](0) = 0$,
the inequality implies  $[v^* (\sum_{i=1}^N u_i)] \equiv 0$ over $(-1/2, 1/2)$, a
contradiction to the supposition that $[v^* (\sum_{i=1}^N u_i)](t) > 0$
for $t > 0$.  Hence we have $[v^* (\sum_{i=1}^N u_i)](t) > 0$
over $(-1, 1)$.

\end{pf}

The next theorem says there is no kink/corner in any length minimizing 
geodesic in $\Tbar$.

\begin{thm}
Every open \weil\ geodesic segment  in $\Tbar$ is entirely contained in a 
single copy of
\teich\ space.
\end{thm}

\begin{rmk}
Given an open geodesic segment, the particular copy of \teich\ space it lies in
may be $\T$ itself, or one component of the frontier $\Tbdry$.  The theorem 
states that the image of a harmonic/energy-minimizing map from an open interval to $\Tbar$
respects the stratified structure of the \weil\ completed \teich\ space $\Tbar$,in the sense that the interior of the geodesic segment 
meeting but a stratum of $\Tbar$
\end{rmk}
  
\begin{pf}
We will first prove that a path which goes through two distinct divisors has a kink, and
therefore it cannot be a \weil\ geodesic.  We will show this by comparing the
lengths of 
two paths: one through the frontier $\Tbdry$, the other through the interior $\T$.

Now let $P_1$ and $P_2$ be two nodes which are obtained by pinching two distinct 
non-intersecting closed curves $c_1$ and $c_2$ on $\Sigma$. Let $z_i, w_i$
with $i=1,2$ be the coordinate systems such that the set $\{ |z_i| < 1, |w_i| < 1, z_i w_i=0 \}$
describes a neighborhood of the node $P_i$ of the surface 
$\Sigma_{c_1 \cup c_2}$.  Recall from the
previous argument that $\{ |z_i| < 1, |w_i| < 1, z_i w_i= t_i \}$ describes 
a fattened
node, and thus $t_i$ gives us a local coordinate of the completed \teich\ space near the 
point $x_0$ representing the nodal surface ${\Sigma}_{c_1 \cup c_2}$.  We assume that the 
coordinate systems $z_i, w_i$ and $t_i$ are the same as the one used 
in the discussion of the degenerating family of hyperbolic metrics. 
Then the hyperbolic lengths $\lambda_i$ 
of the closed geodesic $c_i$
is given by 
\[
\lambda_i = \frac{2 \pi^2}{- \log |t_i|} + 
O \big( \frac{1}{(- \log |t_i|)^4} \big).
\]
as $ (- \log |t_i|)^{-1}$ goes down to zero.
   From now on, we will denote $\frac{2 \pi^2}{- \log |t_1|}$ by $l_1$ and 
$\frac{2 \pi^2}{- \log |t_2|}$ by $l_2$.  Define a functional $L$ defined 
locally in
the neighborhood of $x_0$ in $\Tbar$ by
\[
L = l_1 + l_2.
\]
Note that as the value $\e$ of $L$ goes to zero, the value of $L$ approximates the values
of $\lambda_1 + \lambda_2$.  

  Now we will proceed to calculate \weil\ lengths of two distinct paths near the point $x_0$ representing the nodal surface with two nodes $P_1$ and $P_2$.
The first path $\sigma_1$ describes a deformation of  Riemann surfaces along which the approximate
length functional $L$ remains constant $\e > 0$, i.e. it starts at a Riemann surface $x_1$
with a node $P_2$ and with a closed geodesic $c_1$ whose length is approximately $\e$,
moves through a family of surfaces where the sum of the hyperbolic lengths of the two closed geodesics  
 $c_1$ and $c_2$ are approximately $\e$, and ends at the surface $x_2$ with the node $P_1$ and with
the closed geodesic $c_2$ of length approximately $\e$.  

The second path $\sigma_2$ is chosen to be the path connecting $x_1$ and $x_2$, which goes through a point $x_0$ in ${\T}_{c_1 \cup c_2}$.  In other words,
first pinch off the closed geodesic $c_1$ to a point while keeping the node $P_2$, and then 
secondly fatten the node $P_2$ till it becomes a closed geodesic of hyperbolic length approximately
$\e$ while keeping the node $P_1$.

We need to justify the choice of the second path $\sigma_2$ among all
other paths connecting the two nodal surfaces, which traverse within
the frontier $\Tbdry$.  As $\sigma_2$ is required to go through 
the frontier \teich\ space ${\T}_{c_1 \cup c_2}$, it is necessary
to show that there is no open geodesic segment of $\sigma_2$ lying entirely
in ${\T}_{c_1 \cup c_2}$. Due to the geodesic convexity of the \teich\
spaces ${\T}_{c_1}$, ${\T}_{c_2}$ and ${\T}_{c_1 \cup c_2}$ as well as
the   
Proposition {\textbf 3}, it follows
that $\sigma_2$ intersects with ${\T}_{c_1 \cup c_2}$ at a single point
$x_0$.  For given a point $q$ in ${\T}_{c_1 \cup c_2}$ and point  $p$
in ${\T}_{c_1}$, the proposition says that the interior of the 
length minimizing path connecting them lies
entirely in ${\T}_{c_1}$, hence if $\sigma_2$ has an open geodesic segment
in ${\T}_{c_1 \cup c_2}$, it cannot be locally length minimizing, hence
it cannot be a \weil\ geodesic.    Among all the points available in ${\T}_{c_1 \cup c_2}$ for
$\sigma_2$ to go through, we may assume that the end points $\sigma_2(0)$
and $\sigma_2(1)$ are sufficiently close to ${\T}_{c_2}$ and ${\T}_{c_1}$
respectively so that the  convex functional $f(x) =
d(\sigma_2(0), x) + d(\sigma_2(1), x)$ defined on ${\T}_{c_1 \cup c_2}$
has a minimum $x_0$.  The functional $f$ is convex since
it is a sum of the two convex functionals. The convexity is due to the
fact that each copy of frontier ${\T}_{c}$ is totally geodesically 
imbedded in $\Tbar$.

Let $x_0$ be the single point at which $\sigma_2$ and ${\T}_{c_1 \cup c_2}$
intersect.  Choose $x_0$ to be the origin of the coordinate system
given by 
\[
( l_1, \theta_1, l_2, \theta_2, t_3, ..., t_{3g-3})
\]
where $l_1 = \frac{2 \pi^2}{- \log |t_1|}$, $ l_2 = 
\frac{2 \pi^2}{- \log |t_1|}$, $\theta_1 = \arg t_1$,
and $\theta_2 = \arg t_2$.  
Then suppose 
\[
p = (\e, \theta_1^1, 0, *, \sqrt{\e} t_3^1, \sqrt{\e} t_4^1,..., \sqrt{\e} t_{3g-3}^1)
\]
\[
q = (0, *, \e, \theta_2^2, \sqrt{\e} t_3^2, \sqrt{\e} t_4^2,..., \sqrt{\e} t_{3g-3}^2).
\]
where $*$ indicates that numbers there do not have a meaning, since 
Fenchel-Nielsen deformation is undefined around a node. 
Denote $(t_3^i,...,t_{3g-3}^i)$ by $T_i \in {\bf C}^{3g-8}$.

\begin{lem}
Near the nodal surface $r$, the \weil\ distance function $d$ is approximated
by the distance $d_0$ induced by the model metric 
$dx^2 + du_1^2 + 1/4 (u_1)^6 d \theta_1^2 $ as follows;
\[
|d(p, q) - d_0 (p, q)| = O \Big( \{ u_1 (p) \}^3 \Big)
\]
where $p$ has coordinates $( u_1(p), \theta_1 (p), t_2 (p), ... , 
t_{3g-4}(p)).$
\end{lem}

\begin{rmk}
In the coordinates above, $\theta_1$ varies over the entire ${\bf R}$.
Any two points in the upper half space $\{ u_1 \geq 0 \}$ whose 
$\theta_1$ coordinate differs by an integral multiple of $2 \pi$
represent the same point in the moduli space $M_g$. Furthermore 
a point $( 0, \theta_1, t_2, ... , t_{3g-3})$ 
represents a single nodal surface in $\Tbar$ 
regardless of the value of $\theta_1 \in {\bf R}$.
\end{rmk}

\begin{pf}[of the lemma]
We first choose an arbitrary smooth path $\sigma (s) = (u_1 (s),
\theta_1 (s), t_1 (s),..., t_{3g-3} (s))$ parametrized by arc-length
with respect to the model distance function $d_0$.  We claim then that
the difference between the \weil\ length $L (\sigma)$ of the path $\sigma$
and the $d_0$ length $L_0 (\sigma)$ are a term of size $O(U_1)^3)$ where
$U_1$ is the maximum $u_1$ coordinate $\sigma$ reaches.
\[
\begin{array}{lll}
L ( \sigma ) - L_0 ( \sigma ) & = & \int_0^{L_0} 
\| \sigma' \| - \| \sigma \|_0 ds \\
& = & \int_0^{L_0} \Big[ \{ 1 + O((u_1)^4) \} (u'_1)^2 +  
\{ 1 + O( (u_1)^4 ) \} \sum_i |t'_i|^2  \\
 & & \ \ \ \ + \frac{\{ 1 + O((u_1)^4) \}}{4}
(u_1)^6 (\theta'_1)^2 + \sum_i O((u_1)^3) |u'_1| |t'_i| + 
\sum_i O((u_1)^6) |\theta'_1||t'_i| \Big]^{1/2} \\
 & & \ \ \ \ \ \ \ \  
- \Big[ (u'_1)^2 + \sum_i |t'_i|^2 + \frac{1}{4} (u_1)^6 (\theta'_1)^2 
 \Big]^{1/2}
ds \\
 & = & \int_0^{L_0} \Big[ (u'_1)^2 + \sum_i |t'_i|^2 + \frac{1}{4} (u_1)^6 (\theta'_1)^2 
 + O((u_1)^4) \{ (u'_1)^2 + \sum_i |t'_i|^2 + \frac{1}{4} (u_1)^6 (\theta'_1)^2 \}
 \\
 & & \ \ \ \ + \sum_i O((u_1)^3) |u'_1| |t'_i| + \sum_i O((u_1)^6) |\theta'_1||t'_i| 
 \Big]^{1/2}  \\
 & & \ \ \ \ \ \ \ \ - \Big[ (u'_1)^2 + \sum_i |t'_i|^2 + \frac{1}{4} (u_1)^6 (\theta'_1)^2 
 \Big]^{1/2} ds \\
 & = & \int_0^{L_0} \Big[ O((u_1)^4) +  \sum_i O((u_1)^3) |u'_1| |t'_i| + 
 \sum_i O((u_1)^6) |\theta'_1||t'_i| \Big] ds \\
 & = & O((U_1)^3)
 \end{array}
\]
where we have used the equality $\sqrt{1+x} = 1 + x/2 + o(x)$ as well 
as the fact that $(u'_1)^2 + \sum_i |t'_i|^2 + 
\frac{1}{4}(u_1)^6 (\theta'_1)^2 \equiv 1$
for all $s$, since $\sigma(s)$ is parametrized by arc-length with respect to
$d_0$.  

Now $d(p, q) = L ( \sigma )$ where $\sigma$ is the \weil\ geodesic connecting
$p$ and $q$, while $d_0(p,q) = L_0 (\sigma_0)$ where $\sigma_0$ is 
the $d_0$ length minimizing geodesic.  The estimate above shows that
\[
d_0 (p, q) = L_0 (\sigma_0) = L(\sigma_0) + O \Big( (u_1)^3 \Big) 
\geq d(p,q) + O \Big( (u_1)^3 \Big)
\]
as well as
\[
 d (p, q) = L (\sigma) = L_0 (\sigma) + O \Big( (u_1)^3 \Big) 
\geq d_0 (p,q) + O \Big( (u_1)^3 \Big).
\]
Combining them together we get the statement of the lemma.
\end{pf}
\begin{rmk}
If we are to replace the length with the energy of the path, that is
integrate the square of the differential instead of the first
power, then the corresponding statement is that the difference
between the \weil\ energy and the energy with respect to the model
metric is again within $O \Big( \{ u_1 (p) \}^3 \Big)$ apart.
\end{rmk}

The result above can be easily generalized for the case with two nodes. 
Namely if $d_0$ is the distance function induced by the model metric
\[
ds^2 =  C \frac{dl_1^2}{l_1} + \tilde{C} \frac{dl_2^2}{l_2}  + 
\frac{C}{4} (l_1)^3 d \theta_1^2 + \frac{\tilde{C}}{4} (l_2)^3 d \theta_2^2 +
\sum_{j = 3}^{3g-3} dt_j \otimes d \overline{t_j} 
\]
(In {\bf Lemma 1}, $u_i = \sqrt{l_i}$ was used, which is equivalent
to the above.)
for any smooth path, the \weil\ length of the path is approximated 
by the $d_0$ length of the path with an error of order 
$O \Big( (l_1)^{3/2} \Big) + O \Big( (l_2)^{3/2} \Big)$
and consequently the \weil\ distance function is approximated by $d_0$ as follows.\[
|d(p, q) - d_0 (p, q)| = O \Big( ( \max l_1)^{3/2} \Big) + 
O \Big( (\max l_2)^{3/2} \Big).
\]
where the maximum is taken over the \weil\ geodesic connecting $p$ and
$q$. 

Let ${\sigma_1}$ be a path parametrized as
\[
{\sigma_1}(s) = \Big( (1-s) \e , \theta_1^1, s \e, \theta_2^2, 
\sqrt{\e}\{(1-s)T_1 + s T_2 \} \Big).
\]
with $0 \leq s \leq 1$.

And let ${\sigma_2}$ be a path consisting of two parts
\[
\Big( (1-t) \e, \theta_1^1 , 0, *, \sqrt{\e} (1-t) T_1 \Big)
\]
with $0 \leq t \leq 1$ and
\[
\Big( 0, *, u \e, \theta_2^2, \sqrt{\e} u T_2 \Big)
\]
with $0 \leq u \leq 1$.

Note that each of the two parts of ${\sigma_2}$ is composed of
is a $d_0$ length-minimizing geodesic, which implies that
the $d_0$ length of each path is the same as the $d_0$ distance
between the two end-points. Also note that on $\sigma_1$
the quantity $L = l_1 + l_2$ is held constant $\e$, while
on $\sigma_2$, we have $L \leq \e$.

First calculate the model metric norm of the tangent vector $v(s)$ to the
path ${\sigma_1}(s)$.
\[
v(s) = \Big( -\e, 0, \e, 0, \sqrt{\e}(T_2 - T_1) \Big) = - \e \del_{l_1} + 
\e \del_{l_2} + \sqrt{\e} (t_i^2 - t_i^1) \del_{t_i}
\]
and thus
\[
\begin{array}{lll}
\| v(s) \|_0^2 & = & < v(s), v(s) > \\
          & = & g_{11} (-\e )^2 +  g_{22} (\e)^2 + 2 g_{12} \e(-\e) + 
          | T_2 - T_1 |^2 \e\\
          & = & \frac{C}{l_1} \e^2 +  \frac{\tilde{C}}{l_2} \e^2 + 
          |  T_2 - T_1|^2 \e \\
          & = & \frac{C}{\e (1-s)} \e^2+ \frac{\tilde{C}}{\e s} \e^2 + 
           |  T_2 - T_1 |^2 \e
\end{array}
\]
where $0 \leq s \leq 1$.

Now we integrate the norm of the tangent vector while 
$s$ changes from $0$ to $1/2$.  This corresponds to the
length of the 
first half of the path $\sigma_1$.
\[
\begin{array}{lll}
\int_{0}^{1/2} \|v \|_0 ds    & = & \int_{0}^{1/2} \Big[ \frac{C}{(1-s)} \e + 
      \frac{\tilde{C}}{ s} \e  +  |T_2 - T_2 |^2 \e  \Big]^{1/2}  d s \\
      & = & \sqrt{\e} \int_{0}^{1/2}  \Big[ \frac{C}{1-s} + 
      \frac{\tilde{C}}{s} + |T_2 - T_1 |^2 \Big]^{1/2} d s \\
      & \leq  & \sqrt{\e} \int_{0}^{1/2} \Big[ \frac{C}{1-s}
      + \frac{\tilde{C}}{s} + |T_2|^2 + |T_1|^2 \Big]^{1/2} ds \\
      & < & \sqrt{\e} \int_{0}^{1/2} \Big[ \frac{C}{1-s}
      + |T_1|^2  \Big]^{1/2} ds + \sqrt{\e} \int_{0}^{1/2} 
      \Big[ \frac{\tilde{C}}{s} + |T_2|^2 \Big]^{1/2} ds  \\
      & = &  \sqrt{\e} \int_{1/2}^{1} \Big[ \frac{C}{s}
      + |T_1|^2 \Big]^{1/2} d s + \sqrt{\e} \int_{0}^{1/2} 
      \Big[ \frac{\tilde{C}}{s} + |T_2|^2 \Big]^{1/2} ds 
\end{array}
\]
where the strict inequality comes from the fact that 
$\sqrt{A+B} < \sqrt{A} + \sqrt{B}$ for
$A, B > 0$, and the last equality is due to a change of variable.
By symmetry of the setting with respect to $C$ and $\tilde{C}$, 
we have the following,
\[
\mbox{the $d_0$ length of the second part}    <   \sqrt{\e} \int_{0}^{1/2} 
\Big[ \frac{C}{s}
      + |T_1|^2 \Big]^{1/2} ds + \sqrt{\e} \int_{1/2}^{1} 
      \Big[ \frac{\tilde{C}}{s} + |T_2|^2 \Big]^{1/2} ds.
\]
Hence we have the following estimate for the $d_0$ length $L_0 (\sigma_1)$
of $\sigma_1$.
\[
L_0 ({\sigma_1}) <  \sqrt{\e} \int_{0}^{1} \Big[ \frac{C}{s}
      + |T_1|^2 \Big]^{1/2} ds + \sqrt{\e} \int_{0}^{1} 
      \Big[ \frac{\tilde{C}}{s} + |T_2|^2 \Big]^{1/2} ds. 
\]
On the other hand the $d_0$ length of the path ${\sigma_2}$
is calculated to be 
\[
\begin{array}{lll}
L_0 ({\sigma_2}) & = & \int_0^1 \| v \|_0 dt + \int_0^1 \| v \|_0 du \\
  & = & \int_{0}^{1} \Big[ \frac{C}{(1-t)} \e +
|T_1|^2 \e  \Big]^{1/2} dt +  
\int_{0}^{1} \Big[ \frac{\tilde{C}}{u} \e +
|T_2|^2 \e \Big]^{1/2} du \\
 & = &  \sqrt{\e} \int_{0}^{1} \Big[ \frac{C}{s}
      + |T_1|^2  \Big]^{1/2} d s + \sqrt{\e} \int_{0}^{1} 
      \Big[ \frac{\tilde{C}}{s} + |T_2|^2 \Big]^{1/2} ds
\end{array}
\]
Note here that  by construction, both $L(\sigma_1)$ and $L(\sigma_2)$
are quantities homogeneous in $\e$ of degree $1/2$.
Therefore we have $L({\sigma_1}) < L({\sigma_2})$ for any value of
$\e$.

We claim that any path connecting $p$ and $q$ going through the
frontier \teich\ space ${\T}_{c_1 \cup c_2}$ cannot be length minimizing.
To see this, notice the infimum of the \weil\ distance of such
paths is given by 
\[
L_0({\sigma_2}) + O({\e}^{3/2}) = \sqrt{\e} \Bigg[
 \int_{0}^{1} \Big[ \frac{C}{s}
      + |T_1|^2 \Big]^{1/2} ds + \int_{0}^{1} 
      \Big[ \frac{\tilde{C}}{s} + |T_2|^2 \Big]^{1/2} ds \Bigg] +  O({\e}^{3/2})
\] 
which follows from the length comparison argument above.
While the \weil\ length of the path ${\sigma_1}$ is
bounded strictly from above by the quantity
\[
\sqrt{\e} \Bigg[
 \int_{0}^{1} \Big[ \frac{C}{s}
      + |T_1|^2 \Big]^{1/2} ds + \int_{0}^{1} 
      \Big[ \frac{\tilde{C}}{s} + |T_2|^2 \Big]^{1/2} ds \Bigg] + 
      O({\e}^{3/2}).
\]
This shows that for sufficiently small $\e$ (or equivalently for $p$
and $q$ chosen sufficiently close to the frontier \teich\ space 
${\T}_{c_1 \cup c_2}$,) the \weil\ length of $\sigma_1$ is strictly
less than the \weil\ length of {\it any} path  going entirely
through  the frontier set ${\T}_{c_1} \cup {\T}_{c_2} \cup  
{\T}_{c_1 \cup c_2}$.  

We have so far shown that given two mutually disjoint simple
closed curves $c_1$ and $c_2$, the frontier sets $\overline{{\T}_{c_1}}$
and $\overline{{\T}_{c_2}}$ intersect transversely in the sense
that there is no length minimizing geodesic originating in
${\T}_{c_1}$ ending in ${\T}_{c_2}$ going through ${\T}_{c_1 \cup c_2}$.

To prove the statement of the theorem, this observation needs to
be generalized in the cases where we have two
sets $C_1$ and $C_2$ of simple closed curves, where $C_1 \cup C_2$
represents a collection of mutually disjoint simple closed curves
on $\Sigma$.  

There are two distinct cases, the first being 
when $C_1$ is totally contained in $C_2$.  
Let $p$ be in ${\T}_{C_1}$, $q$ in ${\T}_{C_2}$.
Since ${\T}_{C_2}$ belongs to the  frontier set of ${\T}_{C_1}$, 
by {\bf Lemma 1}, we know that the open geodesic segment connecting
$p$ and $q$ lies entirely in ${\T}_{C_1}$.

The second case is where $C_1 \backslash C_2 \neq \emptyset$ 
and $C_2 \backslash C_1 \neq \emptyset$.  Then we claim that
given $p$ in ${\T}_{C_1}$ and $q$ in ${\T}_{C_2}$, the
open geodesic segment connecting $p$ and $q$ lies in entirely
${\T}_{C_1 \cap C_2}$.  If not, it has to go through  a frontier
\teich\ space ${\T}_{C_1 \cup C_2}$ at a single point $r$.
The \weil\ metric can be approximated by the model metric
\[
\sum_{i=1}^{|C_1|} C_i \Big( du_i^2 + \frac{1}{4}(u_i)^6 d {\theta}_i^2 \Big)     
+ \sum_{j=1}^{|C_2|} C_j \Big( du_j^2 + \frac{1}{4}(u_j)^6 d {\theta}_j^2 \Big)
+ \sum_{k > |C_1| + |C_2|}^{3g-3} d t_k \otimes d \overline{t_k}.
\]
where $u_i = \sqrt{ \frac{2 \pi^2}{ - \log |t_i| }}$ and $\theta_i =
\arg t_i$.  By following the argument for the case when $C_1 = c_1$
and $C_2 = c_2$, it can be checked that the path through $r \in {\T}_{C_1
\cup C_2}$ cannot be length minimizing. 

Finally we have to consider the case 
when we have $C_1$ and $C_2$ be sets of mutually 
disjoint simple closed curves, but $C_1 \cup C_2$ is not
a set of mutually disjoint closed curves, for $p \in {\T}_{C_1}$
and $q \in {\T}_{C_2}$, note a length minimizing geodesic 
lies in ${\T}_{C_1 \cap C_2}$. If not, the open geodesic segment
connecting $p$ and $q$ has to go through at least another
frontier \teich\ space ${\T}_{C_3}$ with $C_3$ strictly containing
$C_1 \cap C_2$.  Then by the argument of the previous paragraph,
there has to be a corner the length minimizing geodesic segment
has to go around, a contradiction.   

\end{pf}

The next theorem had been known, in particular it is a 
consequence of a statement (Theorem 6) which appears in~\cite{Ab}, due to the
fact the the \teich\ distance dominates the \weil\ distance.  The proof is based
on the fact that the Dehn twist can be arbitrarily localized in the
presence of a pinching neck.  The proof is presented here for the sake of
completeness and also to make this idea of localizing the Dehn twist 
explicit in the \weil\ geometric terms.   

It is of particular interest when one studies a local monodromy around a singular fiber (a nodal
surface ${\Sigma}_0$.) ( See for example papers of Matsumoto-Montesinos-Amilibia~\cite{MM}, Earle-Sipe~\cite{ES})

First
define the following functional on $\Tbar$
\[
\delta_{\g} (x) = d(x, \g x).
 \]
Since the \weil\ distance functional $d: \Tbar \times \Tbar$ is strictly
convex in both entries, it follows that $\delta_{\g}$ is a convex functional
as well.  
 
\begin{thm}
Suppose that $\g$ is a Dehn twist around a simple closed geodesic $c$ in $\Sigma$.
Let ${\T}_c$ be the \teich\ space of the surface $\Sigma_0$ obtained by
pinching $c$ of non-singular surface $\Sigma$ to a node.  Then 
the set of points in $\Tbar$ fixed by $\g$ is ${\T}_c$.
\end{thm}

The statement says that the local monodromy is caused
by having a map from the universal cover of a punctured disc
to the \weil\ completed \teich\ space, which is $\rho$ equivariant
where $\rho: {\bf Z} \rightarrow \langle \g \rangle \subset \map.$

\begin{pf}
Suppose \g\ is a Dehn twist around a closed geodesic $c$ on $\Sigma$. Suppose 
$\Sigma_0$ is a Riemann surface with a node $N$ which is obtained by pinching
the closed geodesic $c$.  Then at the node $N$, $\Sigma_0$ has a neighborhood 
isomorphic to $|z| < 1, |w| < 1, zw = 0 $ in ${\bf C}^2$.  

Remove the discs $\{ z : 0 < |z| \leq |t| \}$ and $\{ w : 0 < |w| \leq |t| \}$
from $\Sigma_0$, and then attach $z$ to $t/w$.  Let $A_t = \{ z : |t| < |z| < 1 \}$
and $\alpha$ be the curve $|z| = |w| = |t|^{1/2}$.  Uniformize this 
new non-singular Riemann surface $\Sigma_t$ such that for
a fixed value of $\delta$ over the annular region
$A_{|t|^{1/2}}^{\delta} = \{ z : |t|^{1/2} < |z| < 1-\delta \} \subset A_t$
the conformal factor $\rho$ of the hyperbolic metric $\rho(z) |dz|$ on $\Sigma_t$
satisfy the following uniform estimate;
\[
\frac{1}{C} |z|^{-2} (\log |z|)^{-2} \leq \rho^2 (z) \leq C |z|^{-2} (\log |z|)^{-2}.
\]
for $z$ in $A_{|t|^{1/2}}^{\delta}$, that is $|t|^{1/2} < |z| < 1-\delta$.
This estimate has been improved in the previous section, but 
here it suffices to have the original estimate from~\cite{Ma}.
Now construct a one-parameter family of maps $w_{\theta}$ of $A_t \subset \Sigma_t$ as follows.
\[
w_{\theta} =  
\left\{ \begin{array}{ll}
                       z  &  \mbox{if $|z| < | t |^{3/4}$} \\
                       z \exp \big( i \theta \int_{| t |^{3/4}}^{|z|} 
                       \phi (s) ds \big)
& \mbox{if $| t |^{3/4} \leq  | z | < | t |^{1/4}$} \\
                       z \exp (i \theta ) & \mbox{ if $ | z | > | t |^{1/4}$ }
\end{array}
\right.
\]
where the function $\phi (s)$ is a smooth non-negative function supported on $ | t |^{1/2}
< s < | t |^{1/4} $ with $\int_{| t |^{1/2}}^{| t |^{1/4}} \phi ds = 1$.
Here the number $\theta$ represents the amount of angle the neck has been twisted by.
By differentiating $w_{\theta}$ by $\theta$ at $\theta =0$, we get a vector field 
$\partial_{\theta} w_{\theta}$ on $A_t \subset \Sigma_t$. 
Now take the $\overline{z}$ derivative of this vector field, to obtain the
infinitesimal Beltrami differential $\mu_0 (z)$
\[
\frac{\partial}{\partial \theta} := \partial_{\overline{z}} \partial_{\theta} w_{\theta} = 
i z \phi (| z |) \frac{1}{2}
\frac{z^{1/2}}{{\overline{z}}^{1/2}}
\]

We want to estimate the \weil\ norm of the tangent vector $v_{0}$
induced by the deformation of conformal structure given by this 
infinitesimal Beltrami differential 
$\frac{\partial}{\partial \theta} = \mu_0$.

We recall the geometry of the space of deformations of a given 
hyperbolic metric on a surface. We will use the upper half plane model of the hyperbolic two space ${\bf H}^2$. Then 
the hyperbolic metric has the conformal factor
\[
\rho^2 (z) = \frac{-4}{(z - \overline{z})^2}.
\]
Let $\Gamma$ be the Fuchsian group representing the  Riemann surface $\Sigma$
whose tangent space (to the space of all hyperbolic structures on 
the topological surface) 
the infinitesimal Beltrami differential $\mu_0$ belong to. 
Let $B (\Gamma)$ be the complex Banach space of Beltrami differentials of 
finite $L^{\infty}$ norm which are $\Gamma$ invariant.  A Beltrami
differential $\mu$ is called harmonic if 
$\mu = (z - \overline{z})^2 \overline{\phi}$ for a 
holomorphic quadratic differential.  Let ${\cal B} (\Gamma)$ be the 
subspace of $B (\Gamma)$ consisting of harmonic Beltrami differentials.
Given $\mu$ in $B (\Gamma)$, there is a projection map (see \cite{Wo2}
for example) $P: B (\Gamma)
\rightarrow {\cal B} (\gamma)$ 
\[
P [\mu] = \frac{-3(z - \overline{z})^2}{\pi} \int_{{\bf H}^2} 
\frac{\mu(\zeta)}{(\zeta - \overline{z})^4} d \sigma(\zeta)
\]
where $d \sigma$ is the Euclidean area element.   Now we claim that 
the \weil\ norm of $P[\mu_0]$ is bounded by the $L^2$ norm of $\mu_0$.
Recall that the \weil\ norm $\| \mu_0 \|^2_{WP}$ here is the pairing 
\[
( \mu_0, \phi_{\mu_0} ) = \int_{{\bf H}^2/\Gamma} {\mu}_0 \phi_{\mu} d \sigma (z)
\]
where $\phi_{\mu_0} (z)$ is the holomorphic quadratic differential 
$P[\mu_0] / (z - \overline{z})^2 = P[\mu_0] \rho^2 (z)$, and the $L^2$ norm
$\| \mu_0 \|^2_{L^2}$ is given by
\[
\langle \mu_0, \mu_0 \rangle = \int_{{\bf H}^2/\gamma} |\mu_0|^2 \rho^2(z) 
d \sigma(z).
\]
We denote the term $\mu_0 - P[\mu_0]$ by $\nu_0$.  It was shown by Ahlfors 
(see \cite{Wo2} for example) that $\nu$ belongs to the kernel $N(\Gamma)$ of 
the pairing $B(\gamma) \times QD(\Gamma) \rightarrow {\bf C}$ given by 
$( , )$ defined as above.  To show the claim above, observe 
\[
\begin{array}{lll} 
\int_{{\bf H}^2/\gamma} |{\mu}_0|^2 (z) \rho^2 (z) d \sigma (z) 
 & = & \int_{{\bf H}^2/\gamma} \Big[ (\mu_0 - \nu_0) + \nu_0 \Big] 
 \Big[ \overline{\mu_0 - \nu_0} + \overline{\nu_0} \Big] \rho^2 d \sigma (z) \\ & = & \int_{{\bf H}^2/\gamma} \Big[ | P[\mu_0] |^2 + |\nu_0|^2 + \nu_0 
 \overline{P[\mu_0]} + P[\mu_0] \overline{\nu_0} \Big] \rho^2 (z) d \sigma (z) \\
 & = & \| P[\mu_0] \|_{WP}^2 + \| \nu_0 \|^2_{L^2} \\
  & = & \| v_0 \|_{WP}^2  + \| \nu_0 \|^2_{L^2} \\ 
 \end{array}
 \]
where we have used the fact that $\nu_0$ is perpendicular with $P[\mu_0]$
with respect to the $L^2$ pairing, or equivalently that $\nu_0$ is in the
kernel of the pairing $(\ \ , \phi_{\mu_0} )$. 

Now the \weil\ norm of the tangent vector $v_{0}$ is estimated as follows;
\[
\begin{array}{lll}
\| v_{0} \|^2_{WP} & \leq & \| \mu_0 \|^2_{L^2} \\
 & = & \int_{A_t} |\mu_0|^2 (z) \rho^2 (z) d \sigma (z) \\
 & = & \int_{A_t} \frac{1}{4} |z|^2 \phi^2(z)       \rho^2 (z) d \sigma (z) \\
 & \leq & \int_{A_t} \frac{1}{4} |z|^2 \phi^2(z) C \frac{ d \sigma (z)}{|z|^2 (\log |z|)^2} \\
 & \leq & C \int_0^{2\pi} \int_{| t |^{1/2}}^{| t |^{1/4}} \phi^2(z) 
\frac{|z|dr d\theta}{ (\log |z|)^2} \\
 & \leq & C \frac{|t|^{1/4}}{(\log |t|^{1/4})^2} \int_{| t |^{1/2}}^{| t |^{1/4}}
\phi^2(z) dr \\
 & \leq & C \frac{|t|^{1/4}}{(\log |t|^{1/4})^2} \left( \int_{| t |^{1/2}}^{| t |^{1/4}}
\phi(z) dr \right)^2 \\
 & = & C \frac{|t|^{1/4}}{(\log |t|^{1/4})^2} \rightarrow 0 \mbox{  as $|t| \rightarrow 0$}.
\end{array}
\]

One can also check that the \weil\ norm of the tangent vector  $v_{\theta}$ for any
$\theta \in [0, 2\pi]$ 
also goes down to zero as $|t|$ goes to zero.  

Therefore, the path of a Dehn twist \g\ around a ``fattened'' node of size $|t|$ has 
a \weil\ length  $o(\frac{|t|^{1/4}}{(\log |t|^{1/4})^2}) \rightarrow 0$, 
which in turn implies
that $\delta (\g) \rightarrow 0$ and it is realized on all nodal surfaces with the closed geodesic
$c$ pinched.

We have shown now that $ {\Tbar}_c \subset \{ x \in \Tbar : \g x = x \} $.  We 
will now show the other inclusion.

Suppose $x \notin {\Tbar}_c$ is fixed by $\g$.  Choose a point $y \in {\T}_c$, and
let $\sigma$ be the geodesic connecting $x$ and $y$.  Then note
that $\delta_{\g} (x) = d(x, \g x) = 0 = d (y, \g  y) = \delta_{\g} (y)$.  
Since the number $d(z, \g z)$ is a convex 
functional on $\Tbar$, for any point $w$ on $ \sigma$ we have $d(w, \g w)=0$.
This immediately indicates that $\sigma \subset \Tbdry$, for the Dehn twist
has no fixed point in $\T$.  Since $\sigma$ is a geodesic, it has to lie
entirely in one component ${\T}_c$ of the frontier. However since $x$ is
not in ${\Tbar}_c$, this is not possible, a contradiction.  Hence every
points fixed by $\g$ belongs to ${\Tbar}_c$.

Therefore the statement
\[
{\T}_c = \{ x \in \Tbar : \g x = x \}
\]
follows.
\end{pf}

\begin{cor}
Suppose that $\g$ in $\map$ represents a product of Dehn twists around
a set of mutually disjoint nontrivial closed geodesics $c_1, ... c_n$.
Then the set of fixed points by $\g$ is the frontier \teich\ space 
${\T}_{1...n}$ which represents the collection of Riemann surfaces
obtained by pinching the closed geodesics $c_i$'s.  
\end{cor}

The next theorem had been essentially known (see~\cite{B}) since
the geometry with respect to the \teich\ metric coincides with
that of \weil\ metric in this particular situation. It should be stated 
for the sake of completeness of the picture.

\begin{thm}
Given an element $\g$ of $\map$, there exists a unique point $x$ 
in $\T$ with $\delta_{\g} (x) = 0$ if and only if $\g$ is an element of
finite order.  
\end{thm}
 
The proof follows from the proof of the same statement for \teich\
distance, once one notes that $ d(x, \g x)= 0$
implies the \teich\ distance between $x$ and $\g x$ is also zero.

\section{\weil\ Isometric Action of Mapping Class Groups}


Let us recall Thurston's classification ~\cite{Th} of diffeomorphisms of a Riemann surface.
We will assume the surface is uniformized to have the hyperbolic metric.
An element of $\map = {\cal D}/{\cal D}_0$ is classified as one of the following three types:\\
1) it can be represented by a diffeomorphism of finite order, also called periodic or elliptic;\\
2) it can be represented by a  reducible diffeomorphism, that is, 
the diffeomorphism leaves a tubular neighborhood of 
a collection $C$ of closed geodesics $c_1,... c_n$ invariant;\\
3) it can be represented by a  pseudo-Anosov diffeomorphism (also called irreducible), that is, there is $r>1$ and transverse measured
foliations $F_+, F_{-}$ such that $\g (F_{+}) = r F_{+}$ and 
$\g (F_{-}) = r^{-1} F_{-}$.  In this case the fixed point set of
$\g$ action in ${\cal PMF}( \Sigma)$ (the Thurston boundary of $\T$) 
is precisely $F_+, F_{-}$.\\
\\
As for classification of subgroups, McCarthy and Papadapoulos~\cite{MP} have shown
that  the subgroups of $\map$ is classified into four classes:\\
1) subgroup containing a pair of independent pseudo-Anosov elements
(called {\it sufficiently large subgroups};\\
2) subgroups fixing the pair $\{ F_+ (\g ), F_{-} (\g ) \}$ of fixed points
in ${\cal PMF} (\Sigma)$ for a certain pseudo-Anosov element $\g \in \map$ (such groups
are virtually cyclic);\\
3) finite subgroups;\\
4) infinite subgroups leaving invariant a finite, nonempty, system of
disjoint, non-peripheral, simple closed curves on $\Sigma$
(such subgroups are called reducible.)

We now will relate those classification results with the stratification
structure of the space $\Tbar$.  What one should bear in mind is the correspondence
between various subgroups of a semi-simple Lie group $G$ and totally geodesic
submanifolds of the symmetric space $G/K$ ( where $K$ is the maximal compact
subgroup of $G$) the subgroups stabilize.

\begin{thm} 
Given a reducible element $\g$ of the mapping class group $\map$, leaving a
collection $C$ of mutually disjoint closed geodesics $c_i, i=1, ..., n$ 
invariant, where $n$ is chosen to be maximal.  
Then there is a positive integer $m$ such that ${\g}^m$ stabilizes
the divisor $\cal D$ which represents the collection of nodal surfaces with 
all the $c_i$'s pinched.  
\end{thm}

\begin{rmk} 
Note that  the action of $\g^m$ on each ${\Sigma}_i$  is either finite
or irreducible, for if not, one can introduce an additional node 
which is kept invariant by $\g^m$.
\end{rmk}

\begin{thm}
Given a subgroup $\Gamma$ of the mapping class group $\map$, every element
of which fixes a set $C$ of mutually disjoint closed geodesics $\{ c_i \}$, there
is a subgroup ${\Gamma}' \subset \Gamma$ of a finite index which 
stabilizes the divisor which represents the collection of nodal surfaces
with all the $c_i$'s pinched.  
\end{thm}

Note that the first theorem is a special case of the second theorem, when 
the subgroup is the cyclic group generated by the reducible element $\g$.
We will prove hence the second theorem now.

\begin{pf}
Suppose the action of $\Gamma$ is reducible and  is completely 
reduced by a set of mutually disjoint closed geodesics
$C_1,..., C_r$.  Then   $\Gamma$ acts on
\[
  \Sigma_1 \times ... \times \Sigma_n 
\]
where each $\Sigma_i$ is a component of $\Sigma \backslash \{ C_1 \cup ... \cup
C_r \}$.  Now each element of $\Gamma$ gives a permutation of the connected 
components $\{ \Sigma_1, ... \Sigma_n \}$.  Hence the representation
$\rho: \Gamma \rightarrow \map$ induces a homomorphism
\[
\phi : \Gamma \rightarrow {\cal S}_n
\]
where ${\cal S}_n$ is the symmetric group of $n$ elements.  Let ${\Gamma}'$ be
the kernel of $\phi$.  Then ${\Gamma}'$ is a subgroup of $\Gamma$ of finite 
index, each element of which leaves each punctured surface ${\Sigma}_i$ 
invariant.
This means that regarding the collection of ${\Sigma}_i$'s as a single surface ${\Sigma}_0$
connected by nodes, ${\Gamma}'$ acts on ${\Sigma}_0$ leaving the partition by the nodes
invariant, and therefore the action of ${\Gamma}'$ on $\Tbar$ leaves the divisor$\cal D$ 
representing the nodal surfaces ${\Sigma}_0$ of various conformal structures
invariant.  Note here that we have used the fact that 
the action of $\map$ on $\Tbar$ leaves $\Tbdry$ and $\T$ invariant, consequently that a nodal surface is sent to another nodal surface.  
\end{pf}   

Finally we prove the following statement, which characterizes the action
of pseudo-Anosov element analogous to the isometric action
of a hyperbolic element in ${\rm SL}(2, {\bf R})$.
\begin{thm}
Suppose $\g$ is a pseudo-Anosov element in $\map$ where $\Sigma$ is a surface
possibly with punctures.  Then there exists a $\g$-invariant \weil\ geodesic
in $\T$ of $\Sigma$.  
\end{thm}

\begin{pf}
First we demonstrate, after the argument used in~\cite{B}, 
that it suffices to show that there exists a
point $q$ in $\T$ such that
\[
d(q, \gamma q) = \inf_{x \in \T} d(x, \gamma x).
\]
Suppose we have such a point $q$. Then 
since $\gamma$ is of infinite order, $q$, $\gamma q$ and $\gamma^2 q$
are all distinct points. Let $m_1$ be the mid-point of the 
geodesic segment $\overline{q (\gamma q)}$, $m_2$ that of 
$\overline{(\gamma q) (\gamma^2 q)}$.  Then 
\[
d(q, m_1) = d(m_1, \gamma q) = \frac{1}{2} d(q, \gamma q) 
=  \frac{1}{2} \inf_{x \in \T} d(x, \gamma x).
\]
Similarly 
\[
d(\gamma q, m_2) = d(m_2, \gamma^2 q) = \frac{1}{2} 
d(\gamma q, \gamma^2 q) = 
\frac{1}{2} \inf_{x \in \T} d(x, \gamma x).
\]
and 
\[
\gamma m_1 = m_2.
\]
Then by the triangle inequality,
\[
d(m_1, m_2) \leq d(m_1, \gamma q) + d(\gamma q, m_2) 
= \inf_{x \in \T} d(x, \gamma x).
\]
while we have
\[
d(m_1, m_2) = d(m_1, \gamma m_1) \geq 
\inf_{x \in \T} d(x, \gamma x)
\]
due to the facts that $m_1$, $m_2$ are in $\T$
and that $\overline{m_1 m_2}$ is contained in 
$\T$. (Recall $\T$ is geodesically convex~\cite{Wo1}.) 
Therefore 
\[
d(m_1, m_2) = \inf_{x \in \T} d(x, \gamma x)
\]
which in turn implies that $\gamma q$ is the mid-point of the
geodesic segment $\overline{m_1 m_2}$, and that 
$q$, $\gamma q$ and $\gamma^2 q$ lie on a line $l_1$, where
line here means a harmonic image of ${\bf R}^1$ into $\T$.
The line $l_1$ thus obtained is invariant under the 
isometric action of $\gamma$.

To show that the infimum is achieved at some point $q$ in $\T$, first
let $\{ p_i \}$ be a sequence in ${\T} = {\cal M}_{-1}/{\cal D}_0$ such that 
\[
\lim_{i \rightarrow \infty} d(p_i, \gamma p_i) = \inf_{x \in \T}
d(x, \gamma x). 
\]
Now let $[p_i]$ be the sequence in the moduli space $M_g = 
{\cal M}_{-1}/{\cal D}$. In other words, let $[p_i]$ be the 
image of the projection map $P: {\T}_g \rightarrow M_g$.
Now within the Deligne-Mumford compactification $\overline{M}$
of $M$, find a convergent subsequence $[p_k]$ of $[p_i]$ which
converges to $[p_{\infty}]$ in $\overline{M}$. 

Let $l_i$ be the geodesic loop in $\overline{M}$ based at $p_i$ which lies
in the free homotopy class represented by the pseudo-Anosov
element $\g$.

We state the following technical lemma.

\begin{lem}There is some number $\delta >0$ such that for every $i$,
there is a portion of $l_i$ which lies more than $\delta$ away
from the divisors of the compactified moduli space $\overline{M}$.
\end{lem}

\begin{pf}[of the Lemma]
Suppose not.  Then the sequence of the loops $l_i$ converges to the
divisors as $i$ increases.  In particular, $\lim [p_i] = [p]$ lies in the
set $\cal D$ of the divisors in $\overline{M}$.  Let $q_i$ be a lift 
of $[p_i]$ in $\Tbar$, and suppose $q_i$ converges to $q$, a lift of 
$p$.  Furthermore, as parametrized sets, $l_i(t)$ converges to
a loop $l(t)$ in the boundary set ${\cal D} \subset \overline{M}$. 
Then there exists a collection $C$ of simple closed geodesics 
such that $q$ lies in ${\T}_C$.  Define $g_i$ so that $p_i = g_i q_i$
Then we have 
\[
d(p_i, \g p_i) = d(g_i q_i, \g g_i p_i) = d(q_i, g_i^{-1} q_i g_i q_i).
\]
Since $\g$ is pseudo-Anosov, so is $g_i^{-1} \g g_i$ for each $i$.  We
shall denote $g_i^{-1} \g g_i$ by ${\g}_i$.  

Under the supposition with which we started the proof, for any $\delta >0$  
there exists $N_{\delta}$ so that for any $i > N_{\delta}$, the loop $l_i$ lies entirely in the $\delta$-neighborhood of the divisors.  The geodesic
loop $l_i$ whose initial point as well as end point are $[p_i]$ can
be lifted to $\Tbar$, so that it becomes a \weil\ geodesic segment ${\sigma}_i$
connecting $q_i$ and ${\g}_i q_i$.  Note that $q_i$ converges to $q$ in 
${\T}_C$, while ${\g}_i q_i$
converges to ${\g}_i q$ in ${\T}_{{\g}_i C}$.  The fact that ${\g}_i$ is pseudo-Anosov
implies that for the collection $C$ of mutually disjoint simple closed 
geodesics $c_i$
\[
\Big( \cup_{i = 1}^{|C|} {\Tbar}_{c_i} \Big) \cap 
\g_i \Big(\cup_{i = 1}^{|C|} {\Tbar}_{c_i} \Big) = \emptyset.  
\]

This in turn implies that for each $i > N_{\delta}$, there exists a 
set $\tilde{C}_i$ of mutually disjoint simple closed curves and 
some nonempty subset $\hat{C}_i$ of $C$ which satisfy the following
conditions;
\begin{tabbing}
1) $\hat{C}_i \backslash \tilde{C_i} = \emptyset$ \\
2) $\tilde{C_i} \backslash \hat{C}_i = \emptyset$ \\
3) the set of curves $\hat{C}_i \cup \tilde{C_i}$ can be shrunk to nodes concurrently\\
4) the set of curves $C \cup \tilde{C_i}$ cannot be shrunk to nodes concurrently,
\end{tabbing}
such that the \weil\ geodesic segment ${\sigma}_i$ is contained in
the $\delta$ neighborhood of $\overline{\T}_{\hat{C}_i} \cup \overline{\T}_{\hat{C}_i \cup \tilde{C_i}} 
\cup \overline{{\T}_{{\tilde{C}}_i}}$, where $\overline{\T}_{\hat{C}_i \cup \tilde{C_i}}$ is 
a ``corner'' the geodesic passes through the $\delta$-neighborhood of. 

The condition $\hat{C}_i \backslash \tilde{C_i} = \emptyset$ above says that a point in
${\T}_{\tilde{C_i}}$ represents a hyperbolic surface where at least one of 
the curves represented in $\hat{C}_i$ has a strictly positive
hyperbolic length (a node has been ``fattened.'')  And the condition that
the set of curves $\hat{C}_i \cup \tilde{C_i}$ can be shrunk to nodes concurrently
implies that the frontier sets $\overline{\T}_{\hat{C}_i}$ and $\overline{{\T}_{{\tilde{C}}_i}}$
meet at $\overline{\T}_{\hat{C}_i \cap \tilde{C_i}}$.  The condition 4) is equivalent
to saying that $\overline{{\T}_{C}}$ and $\overline{{\T}_{\tilde{C_i}}}$ are disjoint.
It then follows that for sufficiently small $\delta$, and sufficiently large $i$,
$q_i$ lies in the $\delta$-neighborhood of $\overline{{\T}_{C}}$, though not 
in the $\delta$-neighborhood of $\overline{{\T}_{\tilde{C_i}}}$.  Similarly note that
$\g q_i$ lies in the $\delta$-neighborhood of $\overline{{\T}_{\g C}}$,
though not in the $\delta$-neighborhood of $\overline{{\T}_{\tilde{C_i}}}$

Let $u_k$'s with $1 \leq k \leq |\tilde{C}_i \cup \hat{C}_i| $ be 
the coordinate functions as used in the previous sections
defined near the frontier \teich\ space ${\T}_{\tilde{C}_i \cup \hat{C}_i}$.

Let $v_i$ be the one-dimensional harmonic maps $v_i : [0, 1]
\rightarrow \T$ whose image is the geodesic segment $\sigma_i$.
Define the following pulled-back functional $F_i$
\[
F_i (t) = \sum_{k = 1}^{|\tilde{C}_i \cup \hat{C}_i|} v_i^* u_k.
\]

As shown before, for $u_k$'s sufficiently small (less than some $\e >0$
a constant only depends on the genus $g$) this function $F_i(t)$ defined on $[0,1]$ satisfies
the differential inequality
\begin{equation}
F_i''(t) \leq K_i F_i(t)   \label{eq: concave}
\end{equation}
for $t$ satisfying $F_i(t) < \e$.
Now we claim that the constant $K_i = K(\tilde{C_i} \cup \hat{C}_i)$ in the 
inequality does not
depend on $i$. This is because that $\tilde{C_i} \cup \hat{C}_i$ is determined by the \weil\
geometry of the compactified moduli space $\overline{M_g}$ near the 
intersection/corner $[C]$ of divisors that the loop $l_i$ passes by
for large $i$'s, where $[C]$ is the equivalence class of all the
$\tilde{C_i} \cup \hat{C_i}$'s since all the $\tilde{C_i} \cup \hat{C_i}$'s are conjugates of others
by elements of the mapping class group.  Denote the constant by $K$

Choose $t_0$ so that $\lim_{i \rightarrow \infty} F_i(t_0) = 0$.  Note that
such $t_0$ exists by the construction of $F_i$'s..
Choose $a_i$ and $b_i$ be the numbers so that $[a_i, b_i]$  is the largest
connected interval in $[0, 1]$ containing the number $t_0$ with $F_i(a_i) = F_i(b_i) = \e$.
Reparametrize the domain of $v_i$ via translations and dilations of ${\bf R}^1$ 
so that $F_i(0) = F_i(1)= \e$. Note that the inequality~\ref{eq: concave} is the geodesic
equation in disguise, and that the geodesic equation is invariant under the affine 
change of coordinate of the domain ${\bf R}$.   Define the set $S$ to be
\[
S = \{ t \in [0, 1]: \lim_{i \rightarrow \infty} F_i(t) = 0 \},
\] 
which is nonempty.
Note that for each $t \in S$ the Harnack-type inequality from {\bf Lemma 2} 
\[
\inf_{(a_t, b_t)} F_i \geq L \sup_{(a_t, b_t)}  F_i
\]  
holds for some open interval $(a_t, b_t) \subset [0, 1]$ containing $t$, which in turn
implies that $S$ is open. Recall $L$ in the estimate is independent of $i$.  On
the other hand one can check that $S$ is a closed set
using the facts that $F_i$ are continuous in $t$ and that $F_i$ converges
uniformly to some $F_0(t)$ since the loop $l_i$  converges uniformly to some loop 
in $l_0$ in the boundary $\overline{M_g} \backslash M_g$ of the compactified moduli space $\overline{M_g}$.  

Now we have shown that $S \subset [0, 1]$ is open and closed, which implies $S = [0,1]$, a contradiction 
for $ F_i (0) = F_i (1) 
= \e$ for all $i$.

\end{pf}

Here we modify the $p_i$'s (and $p$ subsequently) so that  
for a sufficiently large $N$ the points
$[p_i]$'s on the loops $l_i$'s with $i > N$ lie outside the 
$\delta$-neighborhood 
of the divisors for the $\delta >0$ from the previous lemma, 
and that $[p_i]$ converges to $[p_{\infty}]$ which
again lies outside the $\delta$-neighborhood of the divisors.  
Let $q_i$ be a lift of 
$[p_i]$ in $\T$, chosen so that $\{q_i\}$ converges to a lift $q$ of
$[p]$. Then let $g_i$ be defined by $p_i = g_i q_i$.  Then we have
\[
d(p_i, \g p_i) = d(g_i q_i, \g g_i q_i) = d(q_i, g_i^{-1} \g g_i q_i) 
\rightarrow \inf_{x \in \Tbar} d(x, \g x)
\]
as $i \rightarrow \infty$.  Denote $ g_i^{-1} \g g_i$ by ${\g}_i$.
Now note that for large $i$'s, $q$ is translated by ${\g}_i$ by
a bounded distance;
\[
\begin{array}{lll}
d(q, {\g}_i q) & \leq & d(q, q_i) + d(q_i, {\g}_i q_i) +
d({\g}_i q_i, {\g}_i q) \\
 & \leq & \e + (\inf d(x, \g x) + \e) + \e \\
 & = & \inf d(x, \g x) + 3 \e
\end{array}
\]
Hence the set of points $\{ {\g}_i q \}$ lies in a ball of radius 
$\Big( \inf d(x, \g x) + 3 \e \Big) $ centered at $q$.  Also recall from 
{\bf Lemma 3} that
the set of points $\{ {\g}_i q \}$ lies outside the $\delta$ neighborhood
$N_{\delta}$ of the frontier sets. 

Now denote the unit tangent vector at $q$ of the geodesic $\sigma_i$ connecting $q$ and ${\g}_i q$  
by $w_i$. The sequence $\{ w_i \}$ has a convergent subsequence on the unit tangent sphere at $q$, which
we denote by $\{ w_i \}$ again.  Let $w_{0}$ be the direction $\{ w_i \}$ converges to.  We now claim
that the geodesic $\sigma_0$ obtained by exponentiating $w_0$ at $q$ lying in the ball 
$B_{(\inf d(x, \g x) + 3 \e)} (q)$ does not hit the frontier sets $\del {\T}$.  

Suppose the contrary. Then $\sigma_0$ hits a frontier \teich\ ${\T}_C$ for some $C \neq \emptyset$.
Define $F_i$ to be the pulled-back function 
\[ 
F_i (x) = \sum_k^{|C|} v_i^* u_k
\]
where $v_i: [0, 1] \rightarrow \T$ is the one dimensional harmonic map whose image is
$\sigma_i$, and $u_k$ with $1 \leq k \leq |C|$ is the coordinate function used in the previous sections
defined near the frontier set ${\T}_C$. Now choose $t_0$ so that
$\lim_{i \rightarrow \infty} F_i (t_0) = 0$.  
For a sufficiently small $\e >0$ as well as sufficiently large $i$'s, let $a_i$ and $b_i$ be the numbers
so that $[ a_i, b_i]$ is the largest connected interval in ${\bf R}$ containing $t_0$ 
with $F_i (a_i) = F_i (b_i) = \e$.  Such $a_i$ and $b_i$ exist because the end points of 
each geodesic segment $\sigma_i$ lies at least $\delta > 0$ distance away from the boundary 
set ${\del} {\T}$.  Note that by the geodesic convexity of $\T$, each
$\sigma_i$ lies in $\T$, and thus $F_i > 0$ for each $i$.    Reparametrize 
$v_i$ via dilations and translations of the domain ${\bf R}$ so that $F_i (0) = F_i (1) = \e$.
 
Let $S$ be the set
\[
S = \{ t : \lim_{i \rightarrow \infty} F_i (t) = 0. \} 
\]
a nonempty subset of $[0, 1]$.
Since $F_i$ satisfies the Harnack-type inequality as noted in the proof of {\bf Lemma 4},
and since $F_i$'s are equicontinuous in $t$, the set $S$ is open and closed in $[0, 1]$, 
and hence $S = [0, 1]$, a contradiction to the fact that $F_i (0) = F_i (1) = \e$ for 
all $i$.  

Now we know that the sequence of the directions $\{ \sigma_i'(0) = w_i \}$ converges to a direction $w_0$
and that the geodesic obtained by exponentiating $w_0$ does not hit the boundary set ${\del} {\T}$
within the ball $B_{(\inf d(x, \g x) + 3 \e)} (q)$.  Then there exists a sufficiently small
neighborhood $N$ of $w_0$ on the unit tangent sphere at $q$, such that  
the exponential map $\exp :   N \times [0, (\inf d(x, \g x) + 3 \e)] \rightarrow \T$ is a diffeormorphism.
Thus it follows that the sequence $\{ {\g}_i q \}$ has a convergent subsequence within the image of the
exponential map above, which is a subset of $\T$.  In particular the convergent subsequence, which
we again denote by
$\{ {\g}_i q \}$, is a Cauchy sequence and for a given ${\e}_0 > 0$, there
exists an integer $N$ such that 
\[
d({\g}_i q, {\g}_j q ) = d(q, {\g}_i^{-1} {\g}_j q ) \leq {\e}_0
\]
for $i, j > N$.

Recall here that the mapping class group $\map$ acts properly discontinuously 
away from the frontier set $\Tbdry$.  That is, for $q$ there exists 
$r_q>0$ such that 
\[
B_{r_q} (q) \cap B_{r_q} (g q) = \emptyset 
\]
for all but finitely many $g$'s in $\map$. We choose ${\e}_0$ above 
to be smaller than the $r_q$, and fix $i$ to be $I > N$. Choose a subsequence
of $\{ {\g}_j \}$ (which we denote also by $\{ {\g}_j \}$) such that
${\g}_I^{-1} {\g}_j = {\g}_0$ where ${\g}_0$ is one of the finite set
of elements in $\map$ which moves $B_{r_q} (q)$ no more than $2 r_q$.

Then for all $j > N$, we have ${\g}_j \equiv {\g}_I {\g}_0$. 
Or equivalently
\[
g_j^{-1} \g g_j = {\g}_I {\g}_0 .
\]
Suppose $j, k > N$.  Then note that 
\[
g_j^{-1} \g g_j = g_k^{-1} \g g_k
\]
and that
\[
(g_k g_j^{-1}) \g = \g (g_k g_j^{-1})
\]
so that $(g_k g_j^{-1})$ commutes with $\g$.  Since $\g$ is pseudo-Anosov,  
it then follows that $(g_k g_j^{-1}) = \g^{n}$ for some $n$.  Therefore
$g_j = {\g}^{n_j} \overline{g}$ for some $\overline{g}$. 
Recall how $p_j$'s were chosen;
\[
d( p_j, \g p_j) \rightarrow \inf_{x \in \Tbar} d(x, \g x)
\]
as $j \rightarrow \infty$. Therefore 
\[
d( p_j, \g p_j) = d(g_j q_j, \g g_j q_j) = 
d({\g}^{n_j} \overline{g} q_j, \g {\g}^{n_j} \overline{g} q_j)
= d(\overline{g} q_j, \g \overline{g} q_j) \rightarrow 
\inf_{x \in \Tbar} d(x, \g x).
\]
Since $\{ q_j \}$ converges to $q$, it follows that $\{ \overline{g} q_j \}$
converges to $ \overline{g} q$ in $\T$.  Therefore we have 
\[
d(q, \g q) = \inf_{x \in \Tbar} d(x, \g x).
\]

\end{pf}

\end{document}